\pdfoutput=1
\RequirePackage{ifpdf}
\ifpdf % We are running pdfTeX in pdf mode
\documentclass[pdftex]{sigma}
\else
\documentclass{sigma}
\fi

\numberwithin{equation}{section}

\newtheorem{Theorem}{Theorem}[section]

\newtheorem{Lemma}[Theorem]{Lemma}
\newtheorem{Proposition}[Theorem]{Proposition}
\newtheorem{Conjecture}[Theorem]{Conjecture}
{\theoremstyle{definition}

\newtheorem{Example}[Theorem]{Example}
\newtheorem{Remark}[Theorem]{Remark}
}

\newcommand{\diag}{\mathop{\rm diag}}

\allowdisplaybreaks

\begin{document}

\renewcommand{\PaperNumber}{028}

\FirstPageHeading

\ShortArticleName{Polynomial Relations for $q$-Characters via the ODE/IM Correspondence}

\ArticleName{Polynomial Relations for $\boldsymbol{q}$-Characters\\ via the ODE/IM Correspondence}

\Author{Juanjuan SUN}

\AuthorNameForHeading{J.~Sun}
\Address{Graduate School of Mathematical Sciences, The University of Tokyo,\\
 Komaba, Tokyo 153-8914, Japan}%
\Email{\href{mailto:sunjuan@ms.u-tokyo.ac.jp}{sunjuan@ms.u-tokyo.ac.jp}}

\ArticleDates{Received January 08, 2012, in f\/inal form May 10, 2012; Published online May 15, 2012}

\Abstract{Let $U_q(\mathfrak{b})$ be the Borel subalgebra of a quantum af\/f\/ine algebra
of type $X^{(1)}_n$ ($X=A,B,C,D$).
Guided by the ODE/IM correspondence in quantum integrable models,
we propose conjectural polynomial relations among the $q$-characters of
certain representations of $U_q(\mathfrak{b})$.}

\Keywords{Borel subalgebra; $q$-character; Baxter's $Q$-operator;
ODE/IM correspondence}

\Classification{81R10; 17B37; 81R50}

%%%%%%%%%%%%%%%%%%%%%%%%%%%%%%%%%%%%%%%%%%%%%%%%%%%%%%%%%%%%%%%%%%%%%%%%%%%%%%%%%%%%%

\section{Introduction}

Let $U_q(\mathfrak{g})$ be a quantum af\/f\/ine algebra of type $X^{(1)}_n$ ($X=A,B,C,D$),
and let $U_q(\mathfrak{b})$ be its Borel subalgebra.
In this paper we shall consider the problem of f\/inding
polynomial relations satisf\/ied by the $q$-characters of
the fundamental modules in the sense of \cite{HJ} and related modules.
This problem is intimately related
with that of functional equations for Baxter's $Q$-operators
in quantum integrable models.
In order to motivate the present study let us review this connection.

In quantum integrable systems, one is interested in
the spectra of a commutative family of transfer matrices.
The latter are constructed from the universal $R$ matrix
of a quantum af\/f\/ine algebra,
by taking the trace of the f\/irst component over some f\/inite-dimensional representation
called `auxiliary space'.
When the auxiliary spaces are the Kirillov--Reshetikhin (KR) modules,
the corresponding transfer matrices satisfy an important family of polynomial identities known as the
$T$-system~\cite{KNS,KS}. (For a recent survey on this topic, see \cite{KNS2011}.)
Subsequently the $T$-system has been formulated and proved~\cite{H, Nak} as identities of
$q$-characters.
It has been shown further~\cite{IIKNS} that
the $T$-system is actually the def\/ining relations of the Grothendieck ring, which is a~polynomial ring~\cite{FR},
 of f\/inite-dimensional modules of quantum af\/f\/ine algebras.
Notice that, from the construction by trace, the $q$-characters and the transfer matrices are both
ring homomorphisms def\/ined on the Grothendieck ring.
Since the $q$-characters are injective~\cite{FM, FR},
identities for $q$-characters imply the same identities for transfer matrices.

Baxter's $Q$-operators were f\/irst introduced in
the study of the $8$-vertex model~\cite{Ba}.
Since then they have been recognized as a key tool in classical and quantum integrable systems,
and there is now a vast literature on this subject.
In the seminal paper~\cite{BLZ1997}, Bazhanov, Lukyanov and Zamolodchikov
revealed that the $Q$-operators can also be obtained from the universal $R$ matrix,
provided the auxiliary space is chosen
to be a (generically inf\/inite-dimensional) representation of the Borel subalgebra.
The work \cite{BLZ1997} for $U_q(\widehat{\mathfrak{sl}}_2)$
has been extended by several authors \cite{Sta, BHK,BaTsu,Kojima,Tsu}
to higher rank and supersymmetric cases.

In view of the results mentioned above,
it is natural to ask whether one can f\/ind polynomial relations, analogous to the $T$-system,
for the $q$-characters of $U_q(\mathfrak{b})$ (see def\/initions in \cite{HJ}) (and hence for the $Q$-operators as well).
The goal of this paper is to propose candidates of such identities.

Our idea is to use
the so-called ODE/IM correspondence which relates
the eigenvalues of $Q$-operators and certain  ordinary linear dif\/ferential equations.
The reader is referred e.g.\ to the review \cite{DDT} on this topic.
In \cite{DDMST}, the correspondence is discussed
for general non-twisted af\/f\/ine Lie algebras using scalar (pseudo-)dif\/ferential operators.
In this paper we reformulate their results
in terms of f\/irst-order systems. The general setting (to be explained in Section~\ref{sec:psi-system}) is as follows.

Let ${}^L\mathfrak{g}$ denote the Langlands dual Lie algebra of $\mathfrak{g}$.
Consider the following ${}^L\mathfrak{g}$-valued f\/irst-order linear dif\/ferential operator
\begin{gather}
\mathcal{L}=\frac{d}{dx}-\frac{\ell}{x}+\sum_{i=1}^n e_i+
\big(x^{Mh^\vee}-E\big) e_0 .
\label{Oper}
\end{gather}
Here $e_0,\ldots,e_n$ are the Chevalley generators of ${}^L\mathfrak{g}$,  $\ell$ is a generic element of
the Cartan subalgebra ${}^L\mathfrak{h}\subset {}^L\mathfrak{g}$, $h^\vee$ is the dual Coxeter number of $\mathfrak{g}$,
and  $M>0$, $E\in\mathbb{C}$ are parameters.
For each fundamental representation $V^{(a)}$ of ${}^L\mathfrak{g}$,
there is a basis
$\{\boldsymbol{\chi}^{(a)}_J(x,E)\}$ of  $V^{(a)}$-valued solutions to the equation $\mathcal{L}\chi=0$
characterized by the behavior
\[
\boldsymbol{\chi}^{(a)}_J(x,E)=\mathbf{u}^{(a)}_J x^{\lambda^{(a)}_J}(1+O(x)),
\qquad x\rightarrow 0,
\]
where $\lambda^{(a)}_J$, $\mathbf{u}^{(a)}_J$
are the eigenvalues of $\ell$ and the corresponding eigenvectors.
Hence this basis is labeled by weight vectors of $V^{(a)}$.
There is also a canonical solution $\boldsymbol{\psi}^{(a)}(x,E)$ which has the fastest decay at $x\to+\infty$.
Due to the special choice~\eqref{Oper} of $\mathcal{L}$, the canonical solutions are shown to
satisfy a system of relations similar to the
Pl{\"u}cker relations, called the $\psi$-system in~\cite{DDMST} (see Subsection~\ref{sec:psisystem} below).
Introduce the connection coef\/f\/icients $Q^{(a)}_J(E)\in \mathbb{C}$ by
\[
\boldsymbol{\psi}^{(a)}(x,E)=\sum_{J} Q^{(a)}_J(E)\boldsymbol{\chi}^{(a)}_J(x,E).
\]
Then the $\psi$-system implies a set of polynomial relations among the connection coef\/f\/i\-cients $Q^{(a)}_J(E)$.

We expect the following to be true:

\begin{enumerate}\itemsep=0pt

\item[(1)] For each connection coef\/f\/icient $Q^{(a)}_J(E)$, there is associated a formal power series $\mathcal{Q}^{(a)}_{J,z}$.
Up to a simple multiplicative factor, the
latter is the $q$-character of an irreducible highest $\ell$-weight module of
$U_q\left(\mathfrak{b}\right)$.
In particular, for the highest or lowest weights of $V^{(a)}$ the corresponding modules are
the fundamental modules \cite{HJ} of $U_q(\mathfrak{b})$.

\item[(2)] With the identif\/ication
\[
Q^{(a)}_J(E) \longleftrightarrow \mathcal{Q}^{(a)}_{J,z} ,\qquad
E\longleftrightarrow z ,
\]
the polynomial relations for $Q^{(a)}_J(E)$ implied by the $\psi$-system
hold true also for $\mathcal{Q}^{(a)}_{J,z}$.
\end{enumerate}

We adopt (1), (2) as our working hypothesis for f\/inding relations among the $q$-characters.
A~few remarks are in order here.
We consider separately the cases related to the spin representation
(i.e., $(\mathfrak{g},a)=(C^{(1)}_n,n), (D^{(1)}_n,n-1), (D^{(1)}_n,n)$) and use the letter
$\mathcal{R}^{(a)}_{\varepsilon,z}$
to denote the counterpart of $\mathcal{Q}^{(a)}_{J,z}$.
We expect that $\mathcal{Q}^{(a)}_{J,z}$ or $\mathcal{R}^{(a)}_{\varepsilon,z}$
corresponds to the $q$-character
of an irreducible highest $\ell$-weight module of $U_q(\mathfrak{b})$
whose highest $\ell$-weight is given by formulas~\eqref{ell-wt}, \eqref{Mez1}--\eqref{Mez2},
 \eqref{Mez3}--\eqref{Mez4}.
(In the case~$C^{(1)}_n$, there are problems with this interpretation,
 as explained in Remark~\ref{Remark-C_n}.)
Also the relations for the $\mathcal{Q}^{(a)}_{J,z}$'s or $\mathcal{R}^{(a)}_{\varepsilon,z}$'s
are not exactly the same as those for $Q^{(a)}_J(E)$'s, but
it is necessary to f\/ine-tune the coef\/f\/icients by some power functions independent of~$z$.
The details will be given in Section~\ref{sec:polyrela} below.
That it is
natural to consider all $\mathcal{Q}^{(a)}_{J,z}$
corresponding to general weights of $V^{(a)}$ is a viewpoint
suggested by the work~\cite{Tsu} for type $A$ algebras.

Now let us come to the content of the present work.
As a f\/irst step, we give explicit candidates for~$\mathcal{Q}^{(1)}_{i,z}$
associated with each weight of the vector representation~$V^{(1)}$.
This is done by taking suitable limits of the known
$q$-characters of KR modules given by tableaux sums.
It is known~\cite{HJ} that  for the highest and lowest weights
this procedure indeed gives the irreducible $q$-characters.
As the next step, we def\/ine the formal series $\mathcal{Q}^{(a)}_{J,z}$
for a general node $a$ of the Dynkin diagram. We def\/ine them by Casorati determinants
whose entries are~$\mathcal{Q}^{(1)}_{i,z}$ with suitable shifts of parameters.
We then give candidates of polynomial relations among the
$\mathcal{Q}^{(a)}_{J,z}$
expected from the $\psi$-system.
In the cases related to spin representations, however,
we do not have explicit candidates the $\mathcal{R}^{(a)}_{\varepsilon,z}$'s
in general,
so we only write down the candidates for the relations.
The resulting polynomial relations are given in Proposition \ref{TypeA-relation}, Theorem \ref{QQ-A_n} (type $A^{(1)}_n$), Conjecture \ref{conj-B_n} (type $B^{(1)}_n$), Conjectures \ref{conj-C_n-1}--\ref{conj-C_n-4} (type $C^{(1)}_n$) and Conjectures \ref{conj-D_n-1}, \ref{conj-D_n-2} (type $D^{(1)}_n$),
respectively.

For the type $A$ algebras,
one can check that $\mathcal{Q}^{(1)}_{i,z}$
(not necessarily the highest or lowest ones) are indeed irreducible $q$-characters
of modules given in~\cite{Kojima}.
The polynomial relations
corresponding to the $\psi$-system
can be summarized as a single identity
\begin{align}
\det\left(\mathcal{Q}^{(1)}_{{\nu},q^{-2\mu+n+2}z}x_{{\nu}}^{-\mu+\nu}\right)^{n+1}_{\mu,\nu=1}
=1 ,
\label{Wronskian}
\end{align}
where $x_\nu=e^{\epsilon_\nu}$, $\{\epsilon_\nu\}_{\nu=1}^{n+1}$ being an orthonormal basis related to the simple roots by
$\alpha_\nu=\epsilon_\nu-\epsilon_{{\nu+1}}$.
This relation (as an identity for $Q$-operators) has been known to experts under the name `Wronskian identity'.
We shall give a direct proof that \eqref{Wronskian} is satisf\/ied for $q$-characters
in Appendix~\ref{app2}.

For the other types of algebras, the situation is less satisfactory.
At the moment we do not know the irreducibility of modules corresponding to the $\mathcal{Q}^{(1)}_{i,z}$
given by our procedure (except those corresponding to the highest or lowest vectors).
For the spin representations of $C^{(1)}_n$ and~$D^{(1)}_n$,
explicit formulas for the $\mathcal{R}^{(a)}_{\varepsilon,z}$'s are missing.
More seriously, we have not been able to prove the proposed identities for $q$-characters
by computational methods.
Instead, we support our working hypothesis by performing the following checks:
\begin{enumerate}\itemsep=0pt
\item $\mathfrak{g}=B^{(1)}_2$, proof by hand,
\item $\mathfrak{g}=B^{(1)}_3$, computer check,
\item $\mathfrak{g}=D^{(1)}_4$, computer check,
\item $\mathfrak{g}=B^{(1)}_n,\,C^{(1)}_n,\, D^{(1)}_n$, proof in the limit to ordinary characters.
\end{enumerate}
The
main results of the present paper consist in
formulating the conjectured relations, and in performing the checks mentioned above.

The text is organized as follows.
In Section \ref{sec:Preliminaries}, we prepare basic def\/initions concerning
the Borel subalgebra $U_q(\mathfrak{b})$
of a quantum af\/f\/ine algebra $U_q(\mathfrak{g})$, and collect necessary facts about their representations.
In Section~\ref{sec:psi-system}, we give an account of the $\psi$-system in the ODE/IM correspondence
and indicate how to derive them using the
formulation by f\/irst-order systems.
We note that in~\cite{DDMST} the $\psi$-system
for algebras other than $A$ type is mentioned as conjectures.
In Section \ref{sec:series},
we introduce the series $\mathcal{Q}^{(1)}_{i,z}$ as limits of the $q$-characters of
KR modules. We also def\/ine these series for the other nodes of the Dynkin diagram.
Section~\ref{sec:polyrela} is devoted to the proposals for polynomial relations.
By comparing with the relations for the connection coef\/f\/icients,
we write down the relations for each
type of algebras $A^{(1)}_n$, $B^{(1)}_n$, $C^{(1)}_n$ and $D^{(1)}_n$.
In Section~\ref{sec:conclusion}, we give a summary of our work.

The text is followed by four appendices.
Appendix~\ref{app1}  gives a list of realizations of the dual Lie algebras ${}^L\mathfrak{g}$.
In Appendix~\ref{app2} we give a proof of the Wronskian identity for type~$A^{(1)}_n$.
In Appendix~\ref{app3}, we prove the identities for type $B^{(1)}_n$ in the limit to
ordinary characters. In Appendix~\ref{app4} the same is done for type~$C^{(1)}_n$ and~$D^{(1)}_n$.

\section{Preliminaries}\label{sec:Preliminaries}

In this section we introduce our notation on quantum af\/f\/ine algebras
and their Borel subalgebras, and collect necessary facts that will be used later.
Throughout this paper, we assume that~$q$ is a nonzero complex number which is not a root of unity.

\subsection{Quantum Borel algebras}

Let $\mathfrak{g}$ be an af\/f\/ine Lie algebra associated with a generalized Cartan matrix
$C=(c_{ij})_{0\leq i,j\leq n}$ of non-twisted type.
Let $D=\textup{diag}(d_0,\ldots,d_n)$ be the unique diagonal matrix such that $DC$ is symmetric and~$d_0=1$.
Set $I=\{1,2,\ldots,n\}$, and let $\mathring{\mathfrak{g}}$ denote
the simple Lie algebra with the Cartan matrix $(c_{ij})_{i,j\in I}$.
Let $\{\alpha_i\}_{i\in I}$, $\{\alpha_i^{\vee}\}_{i\in I}$ and $\{\omega_i\}_{i\in I}$ be the simple roots,
simple coroots and the fundamental weights of $\mathring{\mathfrak{g}}$, respectively.
We set $P=\oplus_{i\in I}\mathbb{Z}\omega_i$,
$Q=\oplus_{i\in I}\mathbb{Z}\alpha_i$.

Set $q_i=q^{d_i}$.
We shall use the standard notation
\[
[k]_i=\frac{q_i^k-q_i^{-k}}{q_i-q_i^{-1}} ,\qquad
[n]_i!=\prod^n_{k=1}[k]_i ,\qquad\left[\begin{matrix}n\\k\end{matrix}\right]_i=\frac{[n]_i!}{[k]_i![n-k]_i!} .
\]

The quantum af\/f\/ine algebra $U_q(\mathfrak{g})$
is the $\mathbb{C}$-algebra def\/ined by generators $E_i$, $F_i$, $K_i^{\pm1}$ ($i=0,\ldots,n$) and the relations
\begin{gather*}
 K_iK_i^{-1}=1=K_i^{-1}K_i ,\qquad K_iK_j=K_jK_i ,\\
 K_iE_jK_i^{-1}=q_i^{c_{ij}}E_j ,\qquad K_iF_jK_i^{-1}=q_i^{-c_{ij}}F_j ,\qquad
 [E_i,F_j]=\delta_{ij}\frac{K_i-K_i^{-1}}{q_i-q_i^{-1}} ,\\
 \sum^{1-c_{ij}}_{r=0}\left[\begin{matrix}1-c_{ij}\\r\end{matrix}\right]_i(-1)^rE_i^rE_jE_i^{1-c_{ij}-r}=0 ,\qquad
i\neq j ,\\
 \sum^{1-c_{ij}}_{r=0}\left[\begin{matrix}1-c_{ij}\\r\end{matrix}\right]_i(-1)^rF_i^rF_jF_i^{1-c_{ij}-r}=0 ,\qquad
i\neq j .
\end{gather*}
We do not write the formulas def\/ining the Hopf structure on $U_q(\mathfrak{g})$ since we are not going to use them.

As is well known~\cite{Be, Dr}, $U_q(\mathfrak{g})$
is isomorphic to the $\mathbb{C}$-algebra with generators $x_{i,r}^{\pm}$
$(i\in I$, $r\in\mathbb{Z})$, $k_i^{\pm1}$ $(i\in I)$, $h_{i,r}$
($i\in I$, $r\in\mathbb{Z}\backslash\{0\}$) and central elements $c^{\pm1/2}$, with the following def\/ining relations
\begin{gather*}
 k_ik_i^{-1}=1=k_i^{-1}k_i ,\qquad
c^{1/2}c^{-1/2}=1=c^{-1/2}c^{1/2} ,\qquad
 k_ik_j=k_jk_i,\qquad k_ih_{j,r}=h_{j,r}k_i ,\\
 k_ix^{\pm}_{j,r}k_i^{-1}=q_i^{\pm c_{ij}}x^{\pm}_{j,r} ,\qquad
 [h_{i,r},x^{\pm}_{j,s}]=\pm\frac{1}{r}[rc_{ij}]_ic^{\mp|r|/2}x^{\pm}_{j,r+s} ,\\
 x^{\pm}_{i,r+1}x^{\pm}_{j,s}-q_i^{\pm
c_{ij}}x^{\pm}_{j,s}x^{\pm}_{i,r+1}=q_i^{\pm
c_{ij}}x^{\pm}_{i,r}x^{\pm}_{j,s+1}-x^{\pm}_{j,s+1}x^{\pm}_{i,r} ,\\
 [x^+_{i,r},x^-_{j,s}]=\delta_{ij}\frac{c^{(r-s)/2}\phi^+_{i,r+s}-c^{-(r-s)/2}\phi^-_{i,r+s}}{q_i-q_i^{-1}} ,\\
 \sum_{\pi\in\mathfrak{S}_{1-c_{ij}}}\sum^{1-c_{ij}}_{k=0}(-1)^k\left[\begin{matrix}1-c_{ij}\\k\end{matrix}\right]_ix^{\pm}_{i,r_{\pi(1)}}\cdots
x^{\pm}_{i,r_{\pi(k)}}x^{\pm}_{j,s}x^{\pm}_{i,r_{\pi(k+1)}}\cdots
x^{\pm}_{i,r_{\pi(1-c_{ij})}}=0,\qquad i\neq j ,
\end{gather*}
for all integers $r_j$, where
$\mathfrak{S}_{m}$ is the symmetric group on $m$ letters, and the $\phi^{\pm}_{i,r}$ are given by
\[
\sum^{\infty}_{r=0}\phi^{\pm}_{i,\pm r}u^{\pm r}=k_i^{\pm1}\exp\left(\pm\big(q_i-q_i^{-1}\big)
\sum^{\infty}_{s=1}h_{i,\pm s}u^{\pm s}\right).
\]

By def\/inition, the Borel subalgebra $U_q(\mathfrak{b})$ is the Hopf
subalgebra of $U_q(\mathfrak{g})$ generated by $E_i$, $K_i^{\pm1}$ ($0\leq i\leq n$).
It is known \cite{Ja} that $U_q(\mathfrak{b})$ is isomorphic to the algebra with
generators $E_i$, $K_i^{\pm1}$ ($0\leq i\leq n$) and the def\/ining relations
\begin{gather*}
 K_iK_j=K_jK_i ,\qquad K_iE_jK_i^{-1}=q_i^{c_{ij}}E_j ,\\
 \sum^{1-c_{ij}}_{r=0}\left[\begin{matrix}1-c_{ij}\\r\end{matrix}\right]_i(-1)^rE_i^rE_jE_i^{1-c_{ij}-r}=0 ,\qquad
i\neq j .
\end{gather*}
The algebra $U_q(\mathfrak{b})$ contains $x^+_{i,m}$, $x^-_{i,r}$, $k_i^{\pm1}$
and $\phi^+_{i,r}$, where $i\in I$, $m\geq0$ and $r>0$.

\subsection[Category $\mathcal{O}$ of $U_q(\mathfrak{b})$]{Category $\boldsymbol{\mathcal{O}}$ of $\boldsymbol{U_q(\mathfrak{b})}$}

In this subsection we recall basics about $U_q(\mathfrak{b})$-modules in category $\mathcal{O}$.
For more details see \cite{HJ}.

Denote by $\mathfrak{t}$ the subalgebra of $U_q(\mathfrak{b})$ generated by $\{k^{\pm 1}_i\}_{i\in I}$.
For a $U_q(\mathfrak{b})$-module $V$ and $\lambda\in P$,
set
\[
V_{\lambda}=\big\{v\in V \mid k_iv= q_i^{\lambda(\alpha_i^\vee)}v, \ \forall\, i\in I\big\} .
\]
When $V_\lambda\neq 0$, it is called the weight space of weight $\lambda$.
A $U_q(\mathfrak{b})$-module $V$ is said to be of type~$1$ if~$V=\oplus_{\lambda\in P}V_\lambda$.

A series of complex numbers
$\Psi=(\Psi_{i,r})_{i\in I,r\in\mathbb{Z}_{\geq0}}$ is called an $\ell$-weight
if $\Psi_{i,0}\neq0$ for all $i\in I$.
We denote by $\mathfrak{t}^{*}_{\ell}$ the set of $\ell$-weights.
For a $U_q(\mathfrak{b})$-module $V$ and $\Psi\in \mathfrak{t}_\ell^*$,
the subspace
\[
V(\Psi)=\big\{v\in V \mid \exists\, p>0 ,\  \forall\, i\in I , \ \forall \, m\ge 0 ,\
(\phi^+_{i,m}-\Psi_{i,m})^pv=0\big\}
\]
is called the $\ell$-weight space of $\ell$-weight $\Psi$.

A $U_q(\mathfrak{b})$-module $V$ is said to be a highest $\ell$-weight module
of highest $\ell$-weight $\Psi\in\mathfrak{t}_\ell^*$
if there exists
a vector $v\in V$ such that $V=U_q(\mathfrak{b})\, v$ and
\[
E_i v=0,\quad i\in I ,\qquad\phi^+_{i,r} v=\Psi_{i,r} v,\quad i\in I, \ r\geq0 .
\]
For each $\Psi\in \mathfrak{t}^{*}_{\ell}$ there exists a
unique simple $U_q(\mathfrak{b})$-module of highest $\ell$-weight $\Psi$.
We denote it by $L(\Psi)$.

A highest $\ell$-weight module is of type $1$ if its highest $\ell$-weight $\Psi$ satisf\/ies
\begin{gather}
\Psi_{i,0}=q_i^{p_i}\quad  \text{for some $p_i\in\mathbb{Z}$, $i\in I$}.
\label{type1hw}
\end{gather}
For any non-zero complex numbers $c_i\in \mathbb{C}^\times$, the map
 $E_i\mapsto E_i$, $K_i\mapsto c_iK_i$ ($i=0,\ldots,n$) gives rise to an  automorphism of $U_q(\mathfrak{b})$.
After twisting by such an automorphism,
any highest $\ell$-weight module can be brought to one satisfying the condition \eqref{type1hw}.
We denote by $\mathfrak{t}^*_{\ell,P}$ the set of $\ell$-weights satisfying \eqref{type1hw}.

Set $D(\lambda)=\lambda-Q^+$, $Q^+=\sum_{i\in I}\mathbb{Z}_{\ge0}\alpha_i$.
A $U_q(\mathfrak{b})$-module $V$ of type $1$ is said to be an object in category $\mathcal{O}$ if
\begin{enumerate}\itemsep=0pt
  \item for all $\lambda\in P$ we have $\dim V_{\lambda}<\infty$,
  \item there exist a f\/inite number of elements $\lambda_1,\ldots,\lambda_s\in P$
such that the weights of $V$ are contained in $\bigcup\limits_{j=1,\ldots,s} D(\lambda_j)$.
\end{enumerate}

In what follows we shall identify $\Psi\in\mathfrak{t}^{*}_{\ell}$ with their generating series,
\[
\Psi=(\Psi_1(u),\ldots,\Psi_n(u)) ,\qquad \Psi_i(u)=\sum_{r\geq0}\Psi_{i,r}u^r .
\]
Simple objects in category $\mathcal{O}$ are classif\/ied by the following theorem.

\begin{Theorem}[\cite{HJ}]
Suppose that $\Psi\in\mathfrak{t}^*_{\ell,P}$.
Then the simple module $L(\Psi)$ is an object in
category $\mathcal{O}$ if and only if $\Psi_i(u)$ is a rational function of $u$ for any $i\in I$.
\end{Theorem}

In particular, for $i\in I$ and $z\in\mathbb{C}^{\times}$,
the simple modules $L^{\pm}_{i,z}=L(\Psi)$  def\/ined by the highest $\ell$-weight
\[
\Psi_j(u)=
\begin{cases}
(1-zu)^{\pm1}, &j=i ,\\
1, &j\neq i ,
\end{cases}
\]
are objects in category $\mathcal{O}$.
These modules are called the fundamental representations~\cite{HJ}.

It is known \cite{Bo,CG} that f\/inite-dimensional simple $U_q(\mathfrak{g})$-modules
remain simple when restricted to $U_q(\mathfrak{b})$.
According to the classif\/ication of the former \cite{CP2,CP3}, the simple module
$L(\Psi)$ is f\/inite-dimensional if its highest $\ell$-weight has the form
\[
\Psi_i(u)=q_i^{{\rm deg}P_i}\frac{P_i(q_i^{-1}u)}{P_i(q_iu)},\qquad  \forall\,  i\in I ,
\]
where $P_i(u)$ is a polynomial such that $P_i(0)=1$.
In the case where
\[
P_j(u)=\begin{cases}
\big(1-q_i^{-m+1} z u\big)\big(1-q_i^{-m+3}z u\big)\cdots \big(1-q_i^{m-1} z u\big), &  j=i  ,
\\
1, &  j\neq i  ,
\end{cases}
\]
with some $i\in I$, $m\in\mathbb{Z}_{>0}$ and $z\in \mathbb{C}^\times$,
the module $L(\Psi)$ is called a Kirillov--Resheti\-khin~(KR) module.
We denote it by $W^{(i)}_{m,z}$.

\subsection[Characters and $q$-characters]{Characters and $\boldsymbol{q}$-characters}\label{def-qcharacter}

We recall the def\/inition of $q$-characters (see \cite{HJ}) and characters of representations of $U_q(\mathfrak{b})$.
Let $\mathbb{Z}^{\mathfrak{t}^*_{\ell,P}}$ denote the set of maps $\mathfrak{t}^*_{\ell,P}\to\mathbb{Z}$.
For $\Psi\in \mathfrak{t}^*_{\ell,P}$, def\/ine
$[\Psi]\in \mathbb{Z}^{\mathfrak{t}^*_{\ell,P}}$ by $[\Psi](\Psi')=\delta_{\Psi,\Psi'}$.
For a~$U_q(\mathfrak{b})$-module $V$  of type $1$ in category $\mathcal{O}$,
its $q$-character $\chi_q(V)$ is def\/ined as an element of~$\mathbb{Z}^{\mathfrak{t}^*_{\ell,P}}$,
\[
\chi_q(V)=\sum_{\Psi\in \mathfrak{t}^*_{\ell,P}}\dim V(\Psi)\cdot [\Psi] .
\]
Similarly let $\mathbb{Z}^{P}$ denote the set of maps $P\to \mathbb{Z}$,  and
def\/ine $e^\lambda\in \mathbb{Z}^{P}$ by $e^{\lambda}(\mu)=\delta_{\lambda,\mu}$.
The ordinary character $\chi(V)$ is an element of $\mathbb{Z}^{P}$,
\[
\chi(V)=\sum_{\lambda\in P}\dim V_\lambda\cdot e^\lambda .
\]
We have a natural map $\varpi:\mathbb{Z}^{\mathfrak{t}^*_{\ell,P}}\rightarrow \mathbb{Z}^P$
which sends  $[\Psi]$ to $e^{\lambda}$ such that  $\Psi_{i,0}=q_i^{\lambda(\alpha_i^{\vee})}$. Under $\varpi$
the $q$-character specializes to the ordinary character,
\[
\chi(V)=\varpi\left(\chi_q(V)\right) .
\]

\section[$\psi$-system]{$\boldsymbol{\psi}$-system}\label{sec:psi-system}

In this section, we reformulate the $\psi$-systems given in \cite{DDMST}.

\subsection[${}^L\mathfrak{g}$-connection]{$\boldsymbol{{}^L\mathfrak{g}}$-connection}
From now on, let $\mathfrak{g}$ be an af\/f\/ine Lie algebra of type $X^{(1)}_n$ ($X=A,B,C,D$), and let
${}^L\mathfrak{g}$ denote its Langlands dual algebra.
Let $h^\vee$  be the dual Coxeter number of $\mathfrak{g}$ (see Table~\ref{table1}).

\begin{table}[hbtp]
\centering\caption{\label{table1}}\vspace{1mm}

\begin{tabular}{|l||*{4}{c|}}
\hline
$\mathfrak{g}$       & $A^{(1)}_n$ & $B^{(1)}_n$ & $C^{(1)}_n$   & \tsep{3pt} $D^{(1)}_n$\bsep{1pt}  \\\hline
${}^L\mathfrak{g}$      & $A^{(1)}_n$ & $A^{(2)}_{2n-1}$ & $D^{(2)}_{n+1}$ & \tsep{3pt} $D^{(1)}_n$ \bsep{1pt} \\\hline
$h^\vee$  & $n+1$     &  $2n-1$         &  $n+1$           & \tsep{3pt} $2n-2$ \bsep{1pt}      \\\hline
$\dim V^{(1)}$  & $n+1$     &  $2n$         &  $2n+2$           & \tsep{3pt} $2n$   \bsep{1pt}    \\\hline
\end{tabular}
\end{table}

Denote by  $e_j$, $f_j$, $h_j$ ($0\le j\le n$) the Chevalley generators of ${}^L\mathfrak{g}$.
We set $e=\sum_{j=1}^ne_j$.
Fix also an element $\ell\in {}^L\mathfrak{h}$ from the Cartan subalgebra of ${}^L\mathfrak{g}$, and
let $\zeta\in \mathbb{C}^\times$.
We consider the following ${}^L\mathfrak{g}$-valued connection (cf.\ e.g.~\cite{FFr09}):
\begin{gather}
 \mathcal{L}=\frac{d}{dx}-\frac{\ell}{x}+e+p(x,E)\zeta e_0 ,
\label{Miura}\\
 p(x,E)=x^{Mh^\vee}-E
\qquad  M>0,\quad  E\in \mathbb{C}.
\label{potential}
\end{gather}
Take an element $h_\rho\in{}^L\mathfrak{h}$
such that $[h_\rho,e_j]=e_j$ ($1\leq j\leq n$) and $[h_\rho,e_0]=-(h^\vee-1)e_0$.
The choice  \eqref{potential}
ensures the following symmetry property for $\mathcal{L}=\mathcal{L}(x,E;\zeta)$:
\begin{gather}
\omega^{k(h_\rho-1)}\mathcal{L}\big(x,E;e^{2\pi i k}\zeta\big)\omega^{-k h_\rho}=
\mathcal{L}\big(\omega^kx,\Omega^kE;\zeta\big) ,
\label{Lsymmetry}
\end{gather}
where $k\in \mathbb{C}$ and
\begin{gather}
\omega=e^{\frac{2\pi i}{h^\vee(M+1)}} ,\qquad \Omega=\omega^{h^\vee M} .
\label{omega}
\end{gather}
On any f\/inite-dimensional ${}^L\mathfrak{g}$-module $V$,
\eqref{Miura} def\/ines a f\/irst-order  system of dif\/ferential equations
$\mathcal{L}(x,E;1)\phi(x,E)=0$.
Quite generally, for a ${}^L\mathfrak{g}$-module $V$,
we denote by $V_k$ the ${}^L\mathfrak{g}$-module obtained by twisting $V$
by the automorphism $e_j\mapsto \exp(2\pi i k\delta_{j,0})e_j$,
$f_j\mapsto \exp(-2\pi i k\delta_{j,0})f_j$.
The operator $\mathcal{L}(x,E;e^{2\pi i k})$ represents the action of $\mathcal{L}(x,E;1)$ on $V_k$.
For a $V$-valued solution $\phi(x,E)$ we set
\begin{gather}
\phi_k(x,E)=\omega^{-kh_\rho}\phi\bigl(\omega^k x, \Omega^{k}E \bigr) .
\label{phik}
\end{gather}
Then the symmetry \eqref{Lsymmetry} implies $\mathcal{L}(x,E;e^{2\pi i k})\phi_k(x,E)=0$.

With  each node $a$ of the Dynkin diagram of ${}^L\mathring{\mathfrak{g}}$ is associated
a fundamental module~$V^{(a)}$ of~${}^L\mathfrak{g}$. We summarize our convention about them and some facts
which will be used later. We leave the proofs to the readers (see Subsection~\ref{subsection5.3} for an example).

\begin{Remark}
In \cite{DDMST}, scalar (pseudo-)dif\/ferential equations are considered.
Using the realization of ${}^L\mathfrak{g}$ given in Appendix \ref{app1} and
rewriting the equation $\mathcal{L}\phi=0$ for the highest component of $\phi$,
one obtains the formulas \cite[(3.18)--(3.21)]{DDMST} (for simplicity we have taken $K=1$ there).
\end{Remark}

The module $V^{(1)}$ is called the vector representation of ${}^L\mathfrak{g}$.
Its explicit realization is given in Appendix \ref{app1}.
If $k$ is an integer, it is obvious that $V_k=V$ for any $V$.
The vector representation
for $\mathfrak{g}=C^{(1)}_n$ has the additional property
\begin{gather}
V^{(1)}_{\frac{1}{2}}\simeq V^{(1)},\qquad  \mathfrak{g}=C^{(1)}_n .
\label{Cn1/2}
\end{gather}

For general $a$, we distinguish the following two cases,
\begin{align*}
{\rm(NS):}\quad (\mathfrak{g},a)&\neq (C^{(1)}_n,n) ,\  (D^{(1)}_n,n-1) ,\  (D^{(1)}_n,n) ,
\\
{\rm(S):}\quad (\mathfrak{g},a)&=(C^{(1)}_n,n) ,\  (D^{(1)}_n,n-1) ,\  (D^{(1)}_n,n) .
\end{align*}
We shall refer to them as the non-spin case and the spin case, respectively.
Here and after we set
\[
t=\begin{cases}
2, &  \mathfrak{g}=C^{(1)}_n ,\\
1, &  \text{otherwise} .
\end{cases}
\]
In the non-spin case, we have
\begin{gather}
 V^{(a)}=\bigwedge^a  V^{(1)}_{\frac{a-1}{2t}}
\qquad\text{for (NS)} .\label{caseA}
\end{gather}
In the case $\mathfrak{g}=B^{(1)}_n$, we have in addition
\begin{gather}
V^{(a)}_{\frac{1}{2}}\simeq V^{(2n-a)},\qquad  \mathfrak{g}=B^{(1)}_n, \quad a=1,\ldots,n ,
\label{Bn1/2}
\end{gather}
where $V^{(a)}$ for $a>n$ stands for the right-hand side of~\eqref{caseA}.

In the case  $(\mathfrak{g},a)= (C^{(1)}_n,n)$,
$V^{(n)}$ is the spin representation of the subalgebra ${}^L\mathring{\mathfrak{g}}=\mathfrak{o}(2n+1)\subseteq D^{(2)}_{n+1}$.
Likewise, in the cases $(\mathfrak{g},a)=(D^{(1)}_n,n-1)$ or $(D^{(1)}_n,n)$,
$V^{(n-1)}$ and $V^{(n)}$ are the spin representations
of the subalgebra ${}^L\mathring{\mathfrak{g}}=\mathfrak{o}(2n)\subseteq D^{(1)}_n$.

Let us consider the solutions at the irregular singularity  $x=\infty$.
It is convenient to use a~gauge transformed form of $\mathcal{L}$,
\begin{align}
x^{Mh_\rho}\mathcal{L}\, x^{-Mh_\rho}=\frac{d}{dx}+\Lambda x^M-\frac{\ell+Mh_\rho}{x}-Ee_0 x^{-M(h_\rho-1)} ,
\label{guage}
\end{align}
where $\Lambda=e+e_0$.
Let $\mu^{(a)}$ be the eigenvalue of $\Lambda$ on $V^{(a)}$ which has the largest real part.
This eigenvalue is multiplicity free, and is given explicitly as follows:
\[
\mu^{(a)}=
\begin{cases}
\displaystyle{
\frac{\sin\frac{\pi a}{t h^\vee}}{\sin\frac{\pi}{t h^\vee}}} &
\text{for (NS)} ,
\vspace{2mm}\\
\displaystyle{\frac{1}{2\sin\frac{\pi}{t h^\vee}}} &
\text{for (S)}.
\end{cases}
\]
Let $\mathbf{u}^{(a)}$ be an eigenvector of $\Lambda$ corresponding to $\mu^{(a)}$.
From the representation \eqref{guage} it follows that
there is a unique $V^{(a)}$-valued
solution $\boldsymbol{\psi}^{(a)}(x,E)$
which satisf\/ies the following in a sector containing the positive real axis $x>0$:
\[
\pmb{\psi}^{(a)}(x,E)=
e^{-\frac{\mu^{(a)}}{M+1}x^{M+1}}x^{-Mh_\rho}\big(\mathbf{u}^{(a)}+o(1)\big)
\quad (x\to\infty) .
\]
We call $\pmb{\psi}^{(a)}(x,E)$ the canonical solution.
In view of the relation~\eqref{caseA} and the formula for $\mu^{(a)}$ given above, we have
\begin{gather}
 \pmb{\psi}^{(a)}=\boldsymbol{\psi}^{(1)}_{-\frac{a-1}{2t}}\wedge \pmb{\psi}^{(1)}_{-\frac{a-3}{2t}}
\wedge\cdots\wedge \pmb{\psi}^{(1)}_{\frac{a-1}{2t}}
\qquad\text{for (NS)} .\label{psi-wedge}
\end{gather}
Here $\pmb{\psi}^{(i)}_k=\omega^{-kh_{\rho}}\pmb{\psi}^{(i)}({\omega^kx,\Omega^kE})$ as def\/ined by \eqref{phik}.

\subsection[$\psi$-system]{$\boldsymbol{\psi}$-system}\label{sec:psisystem}

In this subsection we state the $\psi$-system for the canonical solutions
$\pmb{\psi}^{(a)}(x,E)$ introduced above. Let us consider them case by case.
\medskip

\underline{$\mathfrak{g}=A^{(1)}_n$ (${}^L\mathfrak{g}=A^{(1)}_n$)}.\quad
For $a=1,\ldots,n$, we have the embedding of ${}^L\mathfrak{g}$-modules
\begin{gather}
\iota : \ \bigwedge^2 V^{(a)}_{\frac{1}{2}}\hookrightarrow V^{(a-1)}\otimes V^{(a+1)},
\label{iotaA}
\end{gather}
where $V^{(0)}=V^{(n+1)}=\mathbb{C}$.  The explicit expression of the embedding is given in Subsection~\ref{subsection5.1}.
On the space $V^{(a-1)}\otimes V^{(a+1)}$, the functions $\phi=\iota\bigl(\boldsymbol{\psi}^{(a)}_{-\frac{1}{2}}\wedge \boldsymbol{\psi}^{(a)}_{\frac{1}{2}}\bigr)$
and $\phi=\boldsymbol{\psi}^{(a-1)}\otimes \boldsymbol{\psi}^{(a+1)}$ both satisfy
the equation $\mathcal{L}\phi=0$
and have the behavior
\[
\phi(x,E)=O\left(\exp\left(-\frac{2\mu^{(a)}\cos\frac{\pi}{h^\vee}}{M+1}x^{M+1}\right)\right),
\qquad  x\to\infty
\]
in a sector containing $x>0$. Since such a solution is unique upto a constant multiple,
we conclude that (after adjusting the constant multiple)
\begin{gather}
\iota \bigl(\boldsymbol{\psi}^{(a)}_{-\frac{1}{2}}\wedge \boldsymbol{\psi}^{(a)}_{\frac{1}{2}}\bigr)
=\boldsymbol{\psi}^{(a-1)}\otimes \boldsymbol{\psi}^{(a+1)},
\qquad  a=1,\ldots, n .
\label{psiA-1}
\end{gather}
In particular, denoting by $\{\mathbf{u}_j\}_{j=1}^{n+1}$ the standard basis of $V^{(1)}$
we have
\begin{gather}
\pmb{\psi}^{(n+1)}=\mathbf{u}_1\wedge\cdots\wedge \mathbf{u}_{n+1}.\label{psiA-2}
\end{gather}
We call the relations \eqref{psiA-1}, \eqref{psiA-2} the $\psi$-system for $A^{(1)}_n$.

The $\psi$-system for the other types can be deduced by the same argument,
using the relevant embeddings of representations. We obtain relations of the following form.
\medskip

\underline{$\mathfrak{g}=B^{(1)}_n$ (${}^L\mathfrak{g}=A^{(2)}_{2n-1}$)}.\quad
\begin{gather}
 \iota \bigl(\boldsymbol{\psi}^{(a)}_{-\frac{1}{2}}\wedge \boldsymbol{\psi}^{(a)}_{\frac{1}{2}}\bigr)
=\boldsymbol{\psi}^{(a-1)}\otimes \boldsymbol{\psi}^{(a+1)},\qquad
 a=1,\ldots,n-1 ,\label{psiB-1}
\\
 \iota\bigl(\boldsymbol{\psi}^{(n)}_{-\frac{1}{4}}\wedge \boldsymbol{\psi}^{(n)}_{\frac{1}{4}}\bigr)=
\boldsymbol{\psi}^{(n-1)}_{-\frac{1}{4}}\otimes \boldsymbol{\psi}^{(n-1)}_{\frac{1}{4}} .\label{psiB-2}
\end{gather}
Here $\iota$ stands for the embedding of ${}^L\mathfrak{g}$ \eqref{iotaA} or
\[
\bigwedge^2 V^{(n)}_{\frac{1}{4}}\hookrightarrow V^{(n-1)}_{-\frac{1}{4}}\otimes V^{(n-1)}_{\frac{1}{4}},
\]
which follows from  \eqref{iotaA} and \eqref{Bn1/2}.
\medskip

\underline{$\mathfrak{g}=C^{(1)}_n$ (${}^L\mathfrak{g}=D^{(2)}_{n+1}$)}.\quad
\begin{gather}
 \iota \bigl(\boldsymbol{\psi}^{(a)}_{-\frac{1}{4}}\wedge \boldsymbol{\psi}^{(a)}_{\frac{1}{4}}\bigr)
=\boldsymbol{\psi}^{(a-1)}\otimes \boldsymbol{\psi}^{(a+1)},\qquad
 a=1,\ldots,n-2 ,\label{psiC-1}
\\
 \iota\bigl(\pmb{\phi}^{(n-1)}\bigr)
=\boldsymbol{\psi}^{(n)}_{-\frac{1}{2}}\wedge \boldsymbol{\psi}^{(n)}_{\frac{1}{2}},\label{psiC-2}
\\
 \iota\bigl(\pmb{\phi}^{(n)}\bigr)=
\boldsymbol{\psi}^{(n)}_{-\frac{1}{4}}\otimes \boldsymbol{\psi}^{(n)}_{\frac{1}{4}}.\label{psiC-3}
\end{gather}
Here we have set
\begin{gather*}
\pmb{\phi}^{(n-1)} =\boldsymbol{\psi}^{(1)}_{-\frac{n-2}{4}}\wedge \boldsymbol{\psi}^{(1)}_{-\frac{n-4}{4}}
\wedge\cdots\wedge \boldsymbol{\psi}^{(1)}_{\frac{n-2}{4}} ,\qquad
\pmb{\phi}^{(n)} =\boldsymbol{\psi}^{(1)}_{-\frac{n-1}{4}}\wedge \boldsymbol{\psi}^{(1)}_{-\frac{n-3}{4}}
\wedge\cdots\wedge \boldsymbol{\psi}^{(1)}_{\frac{n-1}{4}} ,
\end{gather*}
and $\iota$ stands for an analog of \eqref{iotaA} or the embeddings
\begin{gather*}
 \bigwedge^{n-1}V^{(1)}_{\frac{n-2}{4}}\simeq V^{(n-1)}\hookrightarrow\bigwedge^2 V^{(n)}_{\frac{1}{2}},
\qquad
 \bigwedge^n V^{(1)}_{\frac{n-1}{4}}\hookrightarrow
V^{(n)}_{-\frac{1}{4}}\otimes V^{(n)}_{\frac{1}{4}},
\end{gather*}
where \eqref{Cn1/2} is taken into account.

\medskip

\underline{$\mathfrak{g}=D^{(1)}_n$ (${}^L\mathfrak{g}=D^{(1)}_n$)}.\quad
\begin{gather}
 \iota \bigl(\boldsymbol{\psi}^{(a)}_{-\frac{1}{2}}\wedge \boldsymbol{\psi}^{(a)}_{\frac{1}{2}}\bigr)
=\boldsymbol{\psi}^{(a-1)}\otimes \boldsymbol{\psi}^{(a+1)},\qquad
 a=1,\ldots,n-3 ,\label{psiD-1}
\\
 \iota\bigl(\pmb{\psi}^{(n-2)}\bigr)=\boldsymbol{\psi}^{(n-1)}_{-\frac{1}{2}}\wedge \boldsymbol{\psi}^{(n-1)}_{\frac{1}{2}} ,
 \qquad \iota\bigl(\pmb{\psi}^{(n-2)}\bigr)=\boldsymbol{\psi}^{(n)}_{-\frac{1}{2}}\wedge \boldsymbol{\psi}^{(n)}_{\frac{1}{2}},\label{psiD-2}
\\
 \iota\bigl(\pmb{\phi}^{(n-1)}\bigr)=\boldsymbol{\psi}^{(n-1)}\otimes \boldsymbol{\psi}^{(n)} .\label{psiD-3}
\end{gather}
We have set
\[
\pmb{\phi}^{(n-1)}=\boldsymbol{\psi}^{(1)}_{-\frac{n-2}{2}}\wedge \boldsymbol{\psi}^{(1)}_{-\frac{n-4}{2}}
\wedge\cdots\wedge \boldsymbol{\psi}^{(1)}_{\frac{n-2}{2}}
\]
and $\iota$ stands for the embedding \eqref{iotaA} or
\begin{gather*}
 \bigwedge^{n-2} V^{(1)}_{\frac{n-3}{2}}\simeq V^{(n-2)}\hookrightarrow \bigwedge^2 V^{(n-1)}_{\frac{1}{2}} ,
 \qquad\bigwedge^{n-2} V^{(1)}_{\frac{n-3}{2}}\simeq V^{(n-2)}\hookrightarrow \bigwedge^2 V^{(n)}_{\frac{1}{2}} ,
\\
 \bigwedge^{n-1} V^{(1)}_{\frac{n-2}{2}}\hookrightarrow V^{(n-1)}\otimes V^{(n)} .
\end{gather*}

\subsection{Connection coef\/f\/icients}

Now we introduce the connection coef\/f\/icients $Q^{(a)}_J(E)$.
First let us consider the vector representation $V^{(1)}$.
We choose $\ell$ generic.
Set $\ell=\textup{diag}(\ell_1,\ldots,\ell_N)$ ($N=\dim V^{(1)}$, see Appendix~\ref{app1}),
and let $\mathbf{u}_j$ ($1\le j\le N$) be the corresponding eigenvector of~$\ell$.
Then there is a unique $V^{(1)}$-valued
solution $\boldsymbol{\chi}^{(1)}_j(x,E)$ characterized by the expansion at the origin,
\[
\boldsymbol{\chi}^{(1)}_j(x,E)=x^{\ell_j}\mathbf{u}_j\left(1+O(x)\right),\qquad  x\to 0  .
\]
From the symmetry \eqref{Lsymmetry} of $\mathcal{L}(x,E)$ we f\/ind that
\begin{gather}
\boldsymbol{\chi}^{(1)}_{j,k}(x,E)=\omega^{k\lambda_j}\boldsymbol{\chi}^{(1)}_j(x,E),\qquad  k\in (1/t)\mathbb{Z}  ,
\label{symchi}
\end{gather}
where we have set
\begin{gather}
\ell-h_\rho=\textup{diag}(\lambda_1,\ldots,\lambda_N) .
\label{lambda-s}
\end{gather}
Def\/ine $Q^{(1)}_j(E)\in\mathbb{C}$ by
\begin{gather}
\boldsymbol{\psi}^{(1)}(x,E)=\sum_{j=1}^N Q^{(1)}_j(E)\boldsymbol{\chi}^{(1)}_j(x,E) .
\label{conn-vec}
\end{gather}

From \eqref{symchi} and \eqref{conn-vec}, we have for $k\in (1/t)\mathbb{Z}$
\[
\boldsymbol{\psi}^{(1)}_k(x,E)=\sum_{j=1}^N Q^{(1)}_{j,k}(E) \omega^{k\lambda_j}\boldsymbol{\chi}^{(1)}_{j}(x,E) ,
\]
where $Q^{(1)}_{j,k}(E)=Q^{(1)}_{j}(\Omega^k E)$.
For a sequence $J=(j_1,\ldots,j_a)$,
introduce the notation
\begin{gather*}
 Q^{(a)}_{J}(E)=\det\left(Q^{(1)}_{j_m,\frac{2l-a-1}{2t}}(E)\cdot
\omega^{(2l-a-1)\lambda_{j_m}/{2t}}\right)_{1\leq l,m\leq a} ,
\\
 \boldsymbol{\chi}^{(a)}_{J}(x,E)=\boldsymbol{\chi}^{(1)}_{j_1,-\frac{a-1}{2t}}(x,E)\wedge
\cdots\wedge \boldsymbol{\chi}^{(1)}_{j_a,-\frac{a-1}{2t}}(x,E) .
\end{gather*}
It follows from \eqref{psi-wedge}
that in the non-spin case
\begin{gather}
 \boldsymbol{\psi}^{(a)}(x,E)=\sum_{J} Q^{(a)}_{J}(E)\,\boldsymbol{\chi}^{(a)}_{J}(x,E) ,
\label{conn-wedge}
\end{gather}
where the sum is taken over all  $J=(j_1,\ldots,j_a)$, $1\le j_1<\cdots<j_a\le N$.

One can similarly def\/ine the connection coef\/f\/icients in the spin case as well.

\section[Series $\mathcal{Q}^{(a)}_{J,z}$, $\mathcal{R}^{(a)}_{\varepsilon,z}$]{Series $\boldsymbol{\mathcal{Q}^{(a)}_{J,z}}$, $\boldsymbol{\mathcal{R}^{(a)}_{\varepsilon,z}}$}\label{sec:series}

It has been shown in \cite{HJ} that the fundamental modules $L^\pm_{i,z}$
of the Borel subalgebra arise as certain limits of the KR modules.
In this section we follow the same procedure to obtain
a~family of formal power series $\mathcal{Q}^{(1)}_{i,z}$ associated with each
weight space of the vector representation of ${}^L\mathfrak{g}$.
Up to simple overall multipliers,
those corresponding to the highest or lowest weights  are
the irreducible $q$-characters $\chi_q(L^{\pm}_{i,z})$.
We expect that in general the $\mathcal{Q}^{(1)}_{i,z}$ are also proportional to
irreducible $q$-characters of $U_q(\mathfrak{b})$. (See however Remark \ref{Remark-C_n}.)

Recall in Subsection \ref{def-qcharacter}, for $\Psi\in\mathfrak{t}^*_{\ell,P}$, the element $[\Psi]\in\mathbb{Z}^{\mathfrak{t}^*_{\ell,P}}$ is def\/ined by $[\Psi](\Psi')=\delta_{\Psi,\Psi'}$. Below we shall use the following elements of  $\mathbb{Z}^{\mathfrak{t}^*_{\ell,P}}$:
\begin{gather}
 \mathcal{Y}_{i,z}=[(1,\ldots,\overset{\text{\tiny $i$-th}}{(1-z u)^{-1}},\ldots,1)] ,
\qquad
e^{\omega_i}=[(1,\ldots,\overset{\text{\tiny $i$-th}}{q_i},\ldots,1)] ,
\label{mcYiz}\\
 Y_{i,z}=[(1,\ldots, \overset{\text{\tiny $i$-th}}{q_i\frac{1-q_i^{-1} zu}{1-q_i zu}},\ldots,1)]
=e^{\omega_i}\mathcal{Y}_{i,{q_iz}}\mathcal{Y}^{-1}_{i,q_i^{-1}z},
\label{Yiz}\\
 A_{i,z}=Y_{i,q_i^{-1}z}Y_{i,q_iz}\prod_{1\leq j\leq n,\atop c_{ji}=-1}Y^{-1}_{j,z}
\prod_{1\leq j\leq n,\atop c_{ji}=-2}Y^{-1}_{j,q_j^{-1}z}Y^{-1}_{j,q_jz}
\prod_{1\leq j\leq n,\atop c_{ji}=-3}Y^{-1}_{j,q_j^{-2}z}Y^{-1}_{j,z}Y^{-1}_{j,q^2_jz}.
\label{Aiz}
\end{gather}
Highest $\ell$-weights are monomials in $\mathcal{Y}^{\pm1}_{i,z}$ and $e^{\pm\omega_i}$.
Abusing the notation, for a monomial $M=[\Psi]$
we shall also write $L(M)$ for $L(\Psi)$.

\subsection{The limiting procedure}

Let us illustrate  on examples the procedure for taking the limit.

\begin{Example}[$\mathfrak{g}=A^{(1)}_2$]
We consider f\/irst the case $\mathfrak{g}=A^{(1)}_2$.
Following \cite{KS,NN}, we write
\begin{gather*}
\fbox{$1$}_{\,z}=Y_{1,z} ,\qquad
\fbox{$2$}_{\,z}=Y^{-1}_{1,q^2z}Y_{2,qz} ,
\qquad
\fbox{$3$}_{\,z}=Y^{-1}_{2,q^3z} .
\end{gather*}
Then the $q$-character  $\chi_q\bigl(W^{(1)}_{m,z}\bigr)$ of the KR module $W^{(1)}_{m,z}$
is presented as a sum over the tableaux
\begin{gather*}
 \sum_{k_1+k_2+k_3=m,\atop k_1,k_2,k_3\geq0}
\overbrace{\begin{tabular}{|l|l|l|}
                  \hline
                  1 & $\cdots$ & 1\\
                  \hline
                \end{tabular}
}^{k_1}\hspace{-0.05cm}\overbrace{\begin{tabular}{l|l|l|}
                  \hline
                  2 & $\cdots$ & 2\\
                  \hline
                \end{tabular}
}^{k_2}\hspace{-0.04cm}\overbrace{\begin{tabular}{l|l|l|}
                  \hline
                  3& $\cdots$ & 3\\
                  \hline
                \end{tabular}
}^{k_3}
\end{gather*}
where the $k$-th box from the right carries the parameter $q^{m+1-2k}z$. This can be rewritten further as
\begin{gather*}
 \prod^m_{j=1}Y_{1,q^{m+1-2j}z}\sum_{k_1+k_2+k_3=m,\atop k_1,k_2,k_3\geq0}\prod^{k_2+k_3}_{j=1}A^{-1}_{1,q^{m+2-2j}z}\prod^{k_3}_{j=1}A^{-1}_{2,q^{m+3-2j}z}\\
\qquad{}= \sum_{k_1+k_2+k_3=m,\atop k_1,k_2,k_3\geq0}e^{(k_1-k_2)\omega_1}e^{(k_2-k_3)\omega_2}\\
\qquad{} \times
\left[\left(\frac{(1-q^{-m}z)(1-q^{m-2k_3+2}z)}{(1-q^{m-2(k_2+k_3)}z)(1-q^{m-2(k_2+k_3)+2}z)},
\frac{(1-q^{m-2(k_3+k_2)+1}z)(1-q^{m+3}z)}{(1-q^{m-2k_3+1}z)(1-q^{m-2k_3+3}z)}\right)\right] .
\end{gather*}
Let us consider the limit $m\rightarrow\infty$.
There are three possibilities to obtain meaningful answers,
\begin{gather*}
k_1\rightarrow\infty, \quad k_2,k_3: \ \text{f\/inite},\quad q^{-k_1}\to0,\\
k_2\rightarrow\infty, \quad k_1,k_3: \ \text{f\/inite},\quad q^{-k_2}\to0,\\
k_3\rightarrow\infty, \quad k_1,k_2: \ \text{f\/inite},\quad q^{-k_3}\to0.
\end{gather*}
Writing $e^{\alpha_i}=x_i/x_{i+1}$ and def\/ining for $i=1,2,3$
\begin{gather*}
\mathcal{Q}^{(1)}_{i,z}=
\prod_{j=1}^{i-1}(1-x_j/x_i)\times
\lim_{k_i\rightarrow\infty\atop q^{-k_i}z\to0}
x_i^{-m}(x_1x_2x_3)^{\frac{m}{3}}\chi_q\left(W^{(1)}_{m,q^{-m}z}\right) ,
\end{gather*}
we obtain the result
\begin{gather*}
\mathcal{Q}^{(1)}_{1,z}
= \mathcal{Y}_{1,z}\sum_{k_2,k_3\geq0}
\prod^{k_2+k_3}_{j=1}A^{-1}_{1,q^{2-2j}z}\prod^{k_3}_{j=1}A^{-1}_{2,q^{3-2j}z} ,
\\
\mathcal{Q}^{(1)}_{2,z}
= \mathcal{Y}^{-1}_{1,q^2z}\mathcal{Y}_{2,qz}\sum_{k_3\geq0}\prod^{k_3}_{j=1}A^{-1}_{2,q^{3-2j}z} ,
\qquad
\mathcal{Q}^{(1)}_{3,z}= \mathcal{Y}^{-1}_{2,q^3z} .
\end{gather*}
Up to simple multipliers, they are the irreducible $q$-characters
\begin{gather*}
 \chi_q\bigl(L^-_{1,z}\bigr)=\mathcal{Q}^{(1)}_{1,z} ,
\qquad
\chi_q\bigl(L(\mathcal{Y}^{-1}_{1,q^2z}\mathcal{Y}_{2,qz})\bigr)=\frac{\mathcal{Q}^{(1)}_{2,z}}{1-e^{-\alpha_1}} ,
\\
 \chi_q\bigl(L^+_{2,q^3z}\bigr)=\frac{\mathcal{Q}^{(1)}_{3,z}}{(1-e^{-\alpha_1-\alpha_2)}(1-e^{-\alpha_2})} .
\end{gather*}
Using the explicit construction of modules \cite{Kojima}, it can be checked
that the second one  is the $q$-character of the simple
module whose highest $\ell$-weight corresponds to the monomial $\mathcal{Y}^{-1}_{1,q^2z}\mathcal{Y}_{2,qz}$.
\end{Example}

We give two more examples.

\begin{Example}[$\mathfrak{g}=B^{(1)}_2$]\label{B12}
\begin{gather*}
\mathcal{Q}^{(1)}_{1,z} =\mathcal{Y}_{1,z}\left\{\sum_{k_2,k_{\bar{2}},k_{\bar{1}}\geq0}\prod^{K_{\bar{1}}}_{l=1}A^{-1}_{1,q^{3-2l}z}
\prod^{K_{\bar{2}}}_{l=1}A^{-1}_{2,q^{2-2l}z}\prod^{K_{\bar{2}}}_{l=1}A^{-1}_{2,q^{3-2l}z}\prod^{K_2}_{l=1}A^{-1}_{1,q^{2-2l}z}\right.\\
 \left.\hphantom{\mathcal{Q}^{(1)}_{1,z} =}{}
 +\sum_{k_2,k_{\bar{2}},k_{\bar{1}}\geq0}\prod^{K_{\bar{1}}}_{l=1}A^{-1}_{1,q^{3-2l}z}
\prod^{K_{\bar{2}}}_{l=1}A^{-1}_{2,q^{2-2l}z}\prod^{K_{\bar{2}}+1}_{l=1}A^{-1}_{2,q^{3-2l}z}\prod^{K_2+1}_{l=1}A^{-1}_{1,q^{2-2l}z}\right\} ,\\
\mathcal{Q}^{(1)}_{2,z} =\mathcal{Y}^{-1}_{1,q^2z}\mathcal{Y}_{2,qz}\left\{\sum_{k_{\bar{2}},k_{\bar{1}}\geq0}\prod^{K_{\bar{1}}}_{l=1}A^{-1}_{1,q^{3-2l}z}
\prod^{K_{\bar{2}}}_{l=1}A^{-1}_{2,q^{2-2l}z}\prod^{K_{\bar{2}}}_{l=1}A^{-1}_{2,q^{3-2l}z}\right.\\
 \left.\hphantom{\mathcal{Q}^{(1)}_{2,z} =}{}
 +\sum_{k_{\bar{2}},k_{\bar{1}}\geq0}\prod^{K_{\bar{1}}}_{l=1}A^{-1}_{1,q^{3-2l}z}
\prod^{K_{\bar{2}}}_{l=1}A^{-1}_{2,q^{2-2l}z}\prod^{K_{\bar{2}}+1}_{l=1}A^{-1}_{2,q^{3-2l}z}\right\} ,\\
\mathcal{Q}^{(1)}_{\bar{2},z} =\mathcal{Y}_{1,qz}\mathcal{Y}^{-1}_{2,q^2z}\sum_{k\geq0}\prod^{k}_{l=1}A^{-1}_{1,q^{3-2l}z} ,\qquad
\mathcal{Q}^{(1)}_{\bar{1},z} =\mathcal{Y}^{-1}_{1,q^3z} .
\end{gather*}
Here we set $K_{\bar1}=k_{\bar1}$,  $K_{\bar2}=k_{\bar2}+k_{\bar1}$
and $K_2=k_2+k_{\bar2}+k_{\bar1}$.
We have
\[
\chi_q(L^-_{1,z})=\mathcal{Q}^{(1)}_{1,z} \qquad
\mbox{and}
 \qquad \chi_q(L^+_{1,q^3z})=\displaystyle\frac{\mathcal{Q}^{(1)}_{\bar{1},z}}{(1-e^{-\alpha_1})(1-e^{-\alpha_1-\alpha_2})(1-e^{-\alpha_1-2\alpha_2})}.
\]
\end{Example}

\begin{Example}[$\mathfrak{g}=C^{(1)}_2$] \label{example-C_2}
\begin{gather*}
\mathcal{Q}^{(1)}_{1,z} =\mathcal{Y}_{1,z}\sum_{k_2,k_0,k_{\bar{2}},k_{\bar{1}}\geq0}\prod^{K_{\bar{1}}}_{l=1}A^{-1}_{1,q^{3-l}z}
\prod^{K_{\bar{2}}}_{l=1}A^{-1}_{2,q^{2-l}z}\prod^{k_0}_{l=1}A^{-1}_{2,q^{2-K_{\bar{2}}-2l}z}\prod^{K_2}_{l=1}A^{-1}_{1,q^{1-l}z} ,\\
\mathcal{Q}^{(1)}_{2,z} =\mathcal{Y}^{-1}_{1,qz}\mathcal{Y}_{2,z}\mathcal{Y}_{2,qz}
\sum_{k_0,k_{\bar{2}},k_{\bar{1}}\geq0}\prod^{K_{\bar{1}}}_{l=1}A^{-1}_{1,q^{3-l}z}
\prod^{K_{\bar{2}}}_{l=1}A^{-1}_{2,q^{2-l}z}\prod^{k_0}_{l=1}A^{-1}_{2,q^{2-K_{\bar{2}}-2l}z} ,\\
\mathcal{Q}^{(1)}_{0,z} =\widehat{\mathcal{Q}}^{(1)}_{0,z}+\widehat{\mathcal{Q}}^{(1)}_{\bar0,z} ,\\
\qquad \widehat{\mathcal{Q}}^{(1)}_{0,z}=\mathcal{Y}_{2,qz}\mathcal{Y}^{-1}_{2,q^2z}
\sum_{k_{\bar{2}},k_{\bar{1}}\geq0,\atop K_{\bar{2}}:{\rm even}}\prod^{K_{\bar{1}}}_{l=1}A^{-1}_{1,q^{3-l}z}
\prod^{\frac{K_{\bar{2}}}{2}}_{l=1}A^{-1}_{2,q^{3-2l}z} ,\\
\qquad \widehat{\mathcal{Q}}^{(1)}_{\bar{0},z}=x_2^{-1}\mathcal{Y}^{-1}_{1,qz}\mathcal{Y}_{1,q^2z}\mathcal{Y}_{2,z}\mathcal{Y}^{-1}_{2,q^3z}
\sum_{k_{\bar{2}},k_{\bar{1}}\geq0, \atop K_{\bar{2}}:{\rm odd}}\prod^{K_{\bar{1}}}_{l=1}A^{-1}_{1,q^{3-l}z}
\prod^{\frac{K_{\bar{2}}-1}{2}}_{l=1}A^{-1}_{2,q^{2-2l}z} ,\\
\mathcal{Q}^{(1)}_{\bar{2},z} =\mathcal{Y}_{1,q^2z}\mathcal{Y}^{-1}_{2,q^2z}\mathcal{Y}^{-1}_{2,q^3z}
\sum_{k_{\bar1}\geq0}\prod^{k_{\bar{1}}}_{l=1}A^{-1}_{1,q^{3-l}z} ,\\
\mathcal{Q}^{(1)}_{\bar{1},z} =\mathcal{Y}^{-1}_{1,q^3z} .
\end{gather*}
Here we set $K_{\bar1}=k_{\bar1}$,  $K_{\bar2}=k_{\bar2}+k_{\bar1}$ and $K_2=k_2+2k_0+k_{\bar2}+k_{\bar1}$.
Note that the weight $0$ has multiplicity $2$.
Correspondingly $\mathcal{Q}^{(1)}_{0,z}$ is a sum of two terms
$\widehat{\mathcal{Q}}^{(1)}_{0,z}$ and $\widehat{\mathcal{Q}}^{(1)}_{\bar0,z}$.
These terms cannot be separated in the process of taking the limit.
Since they have highest $\ell$-weights whose ratio is not a monomial of the $A_{i,z}$'s,
$\mathcal{Q}^{(1)}_{0,z}$ cannot be an irreducible $q$-character.
We have
\[
\chi_q(L^-_{1,z})=\mathcal{Q}^{(1)}_{1,z} \qquad
\mbox{and}\qquad
 \chi_q(L^+_{1,q^3z})=\displaystyle\frac{\mathcal{Q}^{(1)}_{\bar{1},z}}{(1-e^{-\alpha_1})(1-e^{-\alpha_1-\alpha_2})(1-e^{-2\alpha_1-\alpha_2})^2}.
\]
\end{Example}

\subsection[Series $\mathcal{Q}^{(1)}_{i,z}$]{Series $\boldsymbol{\mathcal{Q}^{(1)}_{i,z}}$}\label{section4.2}

In order to discuss the general case, let us prepare some notation.
For $\mathfrak{g}=X^{(1)}_n$ ($X=A,B,C,D$),
introduce a parametrization of $\{\alpha_i,\omega_i\}$ by orthonormal vectors $\{\epsilon_i\}$, and
an index set $\mathcal{J}$ with a~partial ordering $\prec$ as follows.
\begin{align*}
A_n^{(1)}:\qquad &\alpha_i=\epsilon_i-\epsilon_{i+1}\quad (1\leq i\leq n) ,
\\
&\omega_i=\epsilon_1+\cdots+\epsilon_i-\frac{i}{n+1}\sum_{j=1}^{n+1}\epsilon_j\quad (1\leq i\leq n) ,
\\
&\mathcal{J}=\{1,2,\ldots,n+1\} ,\quad 1\prec2\prec\cdots\prec n+1 ,
\\
B_n^{(1)}:\qquad &\alpha_i=\epsilon_i-\epsilon_{i+1}\quad (1\leq i\leq n-1) ,\quad \alpha_n=\epsilon_n ,
\\
&\omega_i=\epsilon_1+\cdots+\epsilon_i\quad (1\leq i\leq n-1) ,\quad
\omega_n=\frac{1}{2}(\epsilon_1+\cdots+\epsilon_n) ,
\\
&\mathcal{J}=\{1,\ldots,n,\bar{n},\ldots,\bar{1}\} ,\quad
1\prec\cdots\prec n\prec\bar{n}\prec\cdots\prec\bar{1} ,
\\
C_n^{(1)}:\qquad &\alpha_i=\frac{1}{\sqrt{2}}(\epsilon_i-\epsilon_{i+1})\quad (1\leq i\leq n-1) ,
\quad
\alpha_n=\sqrt{2}\epsilon_n ,
\\
&\omega_i=
\frac{1}{\sqrt{2}}(\epsilon_1+\cdots+\epsilon_i)\quad (1\leq i\leq n) ,
\\
&\mathcal{J}=\{1,\ldots,n,0,\bar{n},\ldots,\bar{1}\} ,
\quad 1\prec\cdots\prec n\prec0\prec\bar{n}\prec\cdots\prec\bar{1} ,
\\
D_n^{(1)}:\qquad &\alpha_i=\epsilon_i-\epsilon_{i+1}\quad (1\leq i\leq n-1) ,\quad
\alpha_n=\epsilon_{n-1}+\epsilon_{n} ,
\\
&\omega_i=\epsilon_1+\cdots+\epsilon_i\quad(1\leq i\leq n-2) ,\\
&\omega_{n-1}=\frac{1}{2}(\epsilon_1+\cdots+\epsilon_{n-1}-\epsilon_n) ,\quad
\omega_{n}=\frac{1}{2}(\epsilon_1+\cdots+\epsilon_{n-1}+\epsilon_n) ,
\\
&\mathcal{J}=
\{1,\ldots,n,\bar{n},\ldots,\bar{1}\} ,\quad 1\prec\cdots\prec n-1\prec{n\atop\bar{n}}\prec\overline{n-1}
\prec\cdots\prec\bar{1} .
\end{align*}
Def\/ine also $x_i$ and $f_{j,k}$ for $j,k\in\mathcal{J}$ by
\begin{align*}
A_n^{(1)}:\qquad &x_i=e^{\epsilon_i}\,\quad (1\leq i\leq n+1) ,\quad
f_{j,k}=\frac{1}{1-x_k/x_j} ,
\\
B_n^{(1)}:\qquad &x_i=e^{\epsilon_i}=x_{\bar{i}}^{-1} \quad  (1\leq i\leq n) ,
\quad f_{j,k}=\frac{1+\delta_{k,\bar{j}}/x_j}{1-x_k/x_j} ,
\\
C_n^{(1)}:\qquad &x_i=e^{\frac{1}{\sqrt{2}}\epsilon_i}=x_{\bar{i}}^{-1}\quad (1\leq i\leq n), \quad
x_0=1 ,\\
&f_{j,k}=\frac{1}{(1-x_k/x_j)(1+\delta_{k,0}x_k/x_j)} ,
\\
D_n^{(1)}:\qquad&x_i=e^{\epsilon_i}=x_{\bar{i}}^{-1} \quad   (1\leq i\leq n) ,
\quad f_{j,k}=\frac{1-\delta_{k,\bar{j}}x_k/x_j}{1-x_k/x_j} .
\end{align*}

In the following, in the sum of the form $\sum_{k_i,\cdots,k_j}$, unless mentioned explicitly, $k_i,\ldots,k_j$  run over all non-negative integers, and we use the abbreviation $K_l=\sum^{j}_{\mu=l}k_{\mu}$ for $i\prec l\prec j$.
In the case  $\mathfrak{g}=C^{(1)}_n$ and $l\prec0$, we set $K_l=\sum^n_{j=l}k_{j}+2k_0+\sum^n_{j=1}k_{\bar{j}}$.

We give below the formula for $\mathcal{Q}^{(1)}_{i,z}$ for each $i\in\mathcal{J}$.

{\noindent\bf \underline{Case $A_n^{(1)}$}:} For $1\leq i\leq n+1$,
\begin{gather*}
\mathcal{Q}^{(1)}_{i,z}=\Phi_{i,z}\sum_{k_{i+1},\ldots, k_{n+1}}\prod^{n}_{j=i}\prod^{K_{j+1}}_{l=1}A^{-1}_{j,q^{j+1-2l}z},
\end{gather*}
where
\begin{gather}
\Phi_{i,z}=\mathcal{Y}_{i-1,q^iz}^{-1}\mathcal{Y}_{i,q^{i-1}z}.\label{PhiA}
\end{gather}

{\noindent\bf \underline{Case $B_n^{(1)}$}:} For $1\leq i\leq n$,
\begin{gather*}
\mathcal{Q}^{(1)}_{i,z} =\Phi_{i,z} \left\{\sum_{k_{i+1},\ldots,k_n,k_{\bar{n}},\ldots,k_{\bar{1}}}
\prod^{n}_{j=i}\prod^{K_{j+1}}_{l=1}A^{-1}_{j,q^{j+1-2l}z}\prod^n_{j=1}\prod^{K_{\bar{j}}}_{l=1}A^{-1}_{j,q^{2n-j-2l}z}\right.\\
 \left.\hphantom{\mathcal{Q}^{(1)}_{i,z} =}{} +\sum_{k_{i+1},\ldots,k_n,k_{\bar{n}},\ldots,k_{\bar{1}}}
\prod^{n}_{j=i}\prod^{1+K_{j+1}}_{l=1}A^{-1}_{j,q^{j+1-2l}z}\prod^n_{j=1}\prod^{K_{\bar{j}}}_{l=1}A^{-1}_{j,q^{2n-j-2l}z}\right\} ,\\
\mathcal{Q}^{(1)}_{\bar{i},z} =\Phi_{\bar{i},z}
\sum_{k_{\overline{i-1}},\ldots,k_{\bar{1}}}
\prod^{i-1}_{j=1}\prod^{K_{\bar{j}}}_{l=1}A^{-1}_{j,q^{2n-j-2l}z} ,
\end{gather*}
where
\begin{gather}
\Phi_{i,z}=\mathcal{Y}_{i-1,q^iz}^{-1}\mathcal{Y}_{i,q^{i-1}z} ,\qquad
\Phi_{\bar{i},z}=\mathcal{Y}_{i-1,q^{2n-i-1}z}\mathcal{Y}_{i,q^{2n-i}}^{-1} .\label{PhiB}
\end{gather}

{\noindent\bf \underline{Case $C_n^{(1)}$}:}
\begin{gather*}
\mathcal{Q}^{(1)}_{i,z} =\Phi_{i,z}
\sum_{k_{i+1},\ldots,k_n,k_0,k_{\bar{n}},\ldots,k_{\bar{1}}}
\prod^{n-1}_{j=i}\prod^{K_{j+1}}_{l=1}A^{-1}_{j,q^{\frac{1+j-2l}{2}}z}\\
\hphantom{\mathcal{Q}^{(1)}_{i,z} =}{} \times\prod^{n-1}_{j=1}
\prod^{K_{\bar{j}}}_{l=1}A^{-1}_{j,q^{\frac{2n+3-j-2l}{2}}z}\prod^{K_{\bar{n}}}_{l=1}A^{-1}_{n,q^{\frac{n+2-2l}{2}}z}
\prod^{k_0}_{l=1}A^{-1}_{n,q^{\frac{n+2-2(K_{\bar{n}}+2l)}{2}}z},\qquad 1\leq i\leq n ,\\
\mathcal{Q}^{(1)}_{0,z} =\widehat{\mathcal{Q}}^{(1)}_{0,z}+\widehat{\mathcal{Q}}^{(1)}_{\bar0,z} ,
\\
\qquad \widehat{\mathcal{Q}}^{(1)}_{0,z}=\Phi_{0,z}
\sum_{k_{\bar{n}},\ldots,k_{\bar{1}},\atop K_{\bar{n}}:{\rm even}}
\prod^{n-1}_{j=1}\prod^{K_{\bar{j}}}_{l=1}A^{-1}_{j,q^{\frac{2n+3-j}{2}-l}z}
\prod^{\frac{K_{\bar{n}}}{2}}_{l=1}A^{-1}_{n,q^{\frac{n+4}{2}-2l}z} ,\\
\qquad \widehat{\mathcal{Q}}^{(1)}_{\bar{0},z}=x_n^{-1}\Phi_{\bar{0},z}
\sum_{k_{\bar{n}},\ldots,k_{\bar{1}},\atop K_{\bar{n}}:{\rm odd}}
\prod^{n-1}_{j=1}\prod^{K_{\bar{j}}}_{l=1}A^{-1}_{j,q^{\frac{2n+3-j}{2}-l}z}
\prod^{\frac{K_{\bar{n}}-1}{2}}_{l=1}A^{-1}_{n,q^{\frac{n+2}{2}-2l}z} ,\\
\mathcal{Q}^{(1)}_{\bar{i},z} =\Phi_{\bar{i},z}
\sum_{k_{\overline{i-1}},\ldots,k_{\bar{1}}}
\prod^{i-1}_{j=1}\prod^{K_{\bar{j}}}_{l=1}A^{-1}_{j,q^{\frac{2n+3-j}{2}-l}z},\qquad 1\leq i\leq n ,
\end{gather*}
where
\begin{gather}
\Phi_{i,z} =\mathcal{Y}_{i-1,q^{\frac{i}{2}}z}^{-1}\mathcal{Y}_{i,q^{\frac{i-1}{2}}z} , \qquad \Phi_{\bar{i},z}=\mathcal{Y}_{i-1,q^{\frac{2n+2-i}{2}}z}\mathcal{Y}^{-1}_{i,q^{\frac{2n+3-i}{2}}z},\qquad 1\leq i\leq n-1  ,\label{PhiC1}\\
\Phi_{n,z} =\mathcal{Y}_{n-1,q^{\frac{n}{2}}z}^{-1}\mathcal{Y}_{n,q^{\frac{n-2}{2}}z}\mathcal{Y}_{n,q^{\frac{n}{2}}z} ,
\qquad \Phi_{\bar{n},z}=\mathcal{Y}_{n-1,q^{\frac{n+2}{2}}z}\mathcal{Y}^{-1}_{n,q^{\frac{n+2}{2}}z}\mathcal{Y}^{-1}_{n,q^{\frac{n+4}{2}}z} ,\label{PhiC2}\\
\Phi_{0,z} =\mathcal{Y}_{n,q^{\frac{n}{2}}z}\mathcal{Y}^{-1}_{n,q^{\frac{n+2}{2}}z} ,\qquad
\Phi_{\bar{0},z}=\mathcal{Y}_{n-1,q^{\frac{n}{2}}z}^{-1}\mathcal{Y}_{n-1,q^{\frac{n+2}{2}}z}
\mathcal{Y}_{n,q^{\frac{n-2}{2}}z}\mathcal{Y}^{-1}_{n,q^{\frac{n+4}{2}}z} .\label{PhiC3}
\end{gather}

{\noindent\bf \underline{Case $D_n^{(1)}$}:}
\begin{gather*}
\mathcal{Q}^{(1)}_{i,z}= \Phi_{i,z}\left\{\sum_{k_{i+1},\ldots,k_{n},\ldots,k_{\bar{1}}}
\prod^{n-2}_{j=1}\prod^{K_{\bar{j}}}_{l=1}A^{-1}_{j,q^{2n-1-j-2l}z}\prod^{K_{\overline{n-1}}}_{l=1}A^{-1}_{n,q^{n-2l}z}
\prod^{n-1}_{j=i}\prod^{K_{j+1}}_{l=1}A^{-1}_{j,q^{j+1-2l}z}\right.\\
\left. \hphantom{\mathcal{Q}^{(1)}_{i,z}=}{} +\sum_{k_{i+1},\ldots,k_{\bar{n}},\ldots,k_{\bar{1}},\atop k_{\bar{n}}\geq1}
\prod^{n-1}_{j=1}\prod^{K_{\bar{j}}}_{l=1}A^{-1}_{j,q^{2n-1-j-2l}z}\prod^{K_{\bar{n}}}_{l=1}A^{-1}_{n,q^{n-2l}z}
\prod^{n-2}_{j=i}\prod^{K_{j+1}}_{l=1}A^{-1}_{j,q^{j+1-2l}z}\right\},\\
\hphantom{\mathcal{Q}^{(1)}_{i,z}=}{} \quad\quad 1\leq i\leq n-1 ,\\
\mathcal{Q}^{(1)}_{n,z}= \Phi_{n,z}\sum_{k_{\overline{n-1}},\ldots,k_{\bar{1}}}
\prod^{n-2}_{j=1}\prod^{K_{\bar{j}}}_{l=1}A^{-1}_{j,q^{2n-1-j-2l}z}\prod^{K_{\overline{n-1}}}_{l=1}A^{-1}_{n,q^{n-2l}z} ,\\
\mathcal{Q}^{(1)}_{\bar{i},z}= \Phi_{\bar{i},z}\sum_{k_{\overline{i-1}},\ldots,k_{\bar{1}}}
\prod^{i-1}_{j=1}\prod^{K_{\bar{j}}}_{l=1}A^{-1}_{j,q^{2n-1-l-2l}z},\qquad  1\leq i\leq n ,
\end{gather*}
where
\begin{gather}
\Phi_{i,z} =\mathcal{Y}_{i-1,q^iz}^{-1}\mathcal{Y}_{i,q^{i-1}z} ,\qquad \Phi_{\bar{i},z}=\mathcal{Y}_{i-1,q^{2n-2-i}z}\mathcal{Y}_{i,q^{2n-1-i}z}^{-1},\qquad 1\leq
i\leq n-2 ,\label{PhiD1}\\
\Phi_{n-1,z} =\mathcal{Y}_{n-2,q^{n-1}z}^{-1}\mathcal{Y}_{n-1,q^{n-2}z}\mathcal{Y}_{n,q^{n-2}z} ,\qquad \Phi_{\overline{n-1},z}=\mathcal{Y}_{n-2,q^{n-1}z}\mathcal{Y}_{n-1,q^{n}z}^{-1}\mathcal{Y}_{n,q^{n}z}^{-1} ,\label{PhiD2}\\
\Phi_{n,z} =\mathcal{Y}_{n-1,q^{n}z}^{-1}\mathcal{Y}_{n,q^{n-2}z} ,\qquad
\Phi_{\bar{n},z}=\mathcal{Y}_{n-1,q^{n-2}z}\mathcal{Y}_{n,q^nz}^{-1} .\label{PhiD3}
\end{gather}

\begin{Remark}
As in \cite{KS,NN}, it is known that the $q$-character of KR module $W^{(1)}_{m,z}$ for $B^{(1)}_n$ can be written in terms of tableau of elementary boxes $\fbox{a}$, $a\in \{1,\ldots,n,0,\bar{n},\ldots,\bar{1}\}$, while the elementary box $\fbox{0}$ is allowed to appear at most once.
When taking limit of the $q$-character of~$W^{(1)}_{m,z}$, we have only $2n$ kinds of meaningful results, then we use the index set $\{1,\ldots,n,\bar{n},\ldots,\bar{1}\}$. By similar reason we use $\{1,\ldots,n,0,\bar{n},\ldots,\bar{1}\}$ as the index set of $C^{(1)}_n$, which looks like opposite to the usual one.
\end{Remark}

\subsection[Series $\mathcal{Q}^{(a)}_{J,z}$]{Series $\boldsymbol{\mathcal{Q}^{(a)}_{J,z}}$}

In this subsection we give the def\/inition of
$\mathcal{Q}^{(a)}_{J,z}$ where $J\in\mathcal{J}^a$ and
$(\mathfrak{g},a)\neq (C^{(1)}_n,n),(D^{(1)}_n,n-1)$, $(D^{(1)}_n,n)$.
We recall that $q_1=q^{1/2}$ for $C^{(1)}_n$ and $q_1=q$ in the other cases.

We f\/ix our convention about the indices as follows.
For $J=(j_1,\ldots,j_a)\in\mathcal{J}^a$, we denote by~$\underset{\sim}{J}$ the underlying set $\{j_1,\ldots,j_a\}\subseteq\mathcal{J}$.
Set further
\begin{gather}
 \overline{J}=(\overline{j_a},\ldots,\overline{j_1}),\\
 J^*=(i_1,\ldots,i_{b}),\quad
\text{with $\underset{\sim}{J^*}=\mathcal{J}\backslash \underset{\sim}{J}$ \ and \ $i_1\prec\cdots\prec i_b$} .
\end{gather}
We say $J$ is increasing if $j_1\prec\cdots\prec j_a$.

For an element $J=(j_1,\ldots,j_a)\in \mathcal{J}^a$, we def\/ine $\mathcal{Q}^{(a)}_{J,z}$ by
\begin{gather}
\mathcal{Q}^{(a)}_{J,z}
 =\det\left(\mathcal{Q}^{(1)}_{j_{\nu},q_1^{-a-1+2\mu}z}x_{j_{\nu}}^{\mu-1}\right)^a_{\mu,\nu=1}\nonumber\\
\hphantom{\mathcal{Q}^{(a)}_{J,z}}{}  =\begin{vmatrix}
\mathcal{Q}^{(1)}_{j_{1},q_1^{1-a}z}x_{j_1}^{0} &\mathcal{Q}^{(1)}_{j_{2},q_1^{1-a}z}x_{j_2}^{0}&\cdots&
\mathcal{Q}^{(1)}_{j_{a},q_1^{1-a}z}x_{j_a}^{0}\\
\mathcal{Q}^{(1)}_{j_1,q_1^{3-a}z}x_{j_1}&\mathcal{Q}^{(1)}_{j_2,q_1^{3-a}z}x_{j_2}&\cdots&\mathcal{Q}^{(1)}_{j_a,q_1^{3-a}z}x_{j_a}\\
\vdots&\vdots&&\vdots\\
\mathcal{Q}^{(1)}_{j_1,q_1^{a-1}z}x_{j_1}^{a-1}&\mathcal{Q}^{(1)}_{j_2,q_1^{a-1}z}x_{j_2}^{a-1}&
\cdots&\mathcal{Q}^{(1)}_{j_a,q_1^{a-1}z}x_{j_a}^{a-1}
\end{vmatrix}.\label{determint-def}
\end{gather}
We set $\mathcal{Q}^{(0)}_{\varnothing,z}=1$.

Note that if the entries of $J\in\mathcal{J}^a$ are not distinct
(i.e.\ $\sharp\underset{\sim}{J}<a$), then $\mathcal{Q}^{(a)}_{J,z}=0$.

Def\/ine
\begin{gather*}
\bar{\mathcal{Q}}^{(a)}_{J}=\varpi\big(\mathcal{Q}^{(a)}_{J,z}\big).
\end{gather*}
Hence by the def\/inition of $\mathcal{Q}^{(a)}_{J,z}$, it is easy to get
\begin{gather}
\bar{\mathcal{Q}}^{(a)}_{J}
=\prod_{i\in J,j\in J^*,\atop i\prec j}
f_{i,j}\prod_{i,j\in J,\atop i\prec j}\left((x_j-x_i)f_{i,j}\right).
\label{Character-Formula}
\end{gather}

For an increasing element $J=(j_1,\ldots,j_a)\in \mathcal{J}^a$, set
\begin{gather}
\Phi_{J,z}=\prod^a_{k=1}\Phi_{j_k,q_1^{a+1-2k}z} ,\label{ell-wt}
\end{gather}
where $\Phi_{j,z}$ are def\/ined in \eqref{PhiA}--\eqref{PhiD3}
with $\mathcal{Y}_{i,z}$ being given in \eqref{mcYiz}.
Let $L(\Phi_{J,z})$
be the unique irreducible $U_q(\mathfrak{b})$-module
with the highest $\ell$-weight $\Phi_{J,z}$.
Except in the case $\mathfrak{g}=C^{(1)}_n$,
we expect that the $q$-character of $L(\Phi_{J,z})$ is given by
\begin{gather*}
 \chi_q(L(\Phi_{J,z}))=\prod\limits_{i\in J,j\in J^*\atop i\succ j}f_{j,i}
\prod\limits_{i,j\in J,\atop i\prec j}\big({-}x_i^{-1}\big)\cdot \mathcal{Q}^{(a)}_{J,z} .
\end{gather*}
In particular, from \eqref{Character-Formula} we expect that
\begin{gather*}
\chi \left(L(\Phi_{J,z})\right)=\prod_{(i,j)\in(J\times J^*)\cup(J^*\times J),\atop i\prec j}f_{i,j}
\times
\begin{cases}
1,& \mathfrak{g}=A^{(1)}_n  ,\\
\prod\limits_{i,\bar{i}\in J,\atop i\prec \bar{i}}\big(1+x_i^{-1}\big), & \mathfrak{g}=B^{(1)}_n  ,\\
\prod\limits_{i,\bar{i}\in J,\atop i\prec \bar{i}}\big(1-x_i^{-2}\big), & \mathfrak{g}=D^{(1)}_n  .
\end{cases}
\end{gather*}

In the case of $\mathfrak{g}=C^{(1)}_n$,
we expect the same to be true for $\mathcal{Q}^{(a)}_{J,z}$
if $j_1,\dots,j_a\neq 0$.
As noted already in Example \ref{example-C_2},
the series $\mathcal{Q}^{(1)}_{0,z}$ does not correspond to an irreducible $q$-character.
We expect rather that $\widehat{\mathcal{Q}}^{(1)}_{0,z}$, $\widehat{\mathcal{Q}}^{(1)}_{\bar0,z}$
correspond to irreducible $q$-characters.
See also Remark \ref{Remark-C_n}.

\subsection[Series $\mathcal{R}^{(n)}_{\varepsilon,z}$ for the spin node: case $C^{(1)}_n$]{Series $\boldsymbol{\mathcal{R}^{(n)}_{\varepsilon,z}}$ for the spin node: case $\boldsymbol{C^{(1)}_n}$}

In this subsection, we introduce another set of series $\mathcal{R}^{(n)}_{\varepsilon,z}$
for the spin node.
Unlike the series~$\mathcal{Q}^{(a)}_{J,z}$,
in general we do not know the explicit formulas for them.
We def\/ine $\mathcal{R}^{(n)}_{\varepsilon,z}$ by the $q$-characters of the irreducible
$U_q(\mathfrak{b})$ module $L(M_{\varepsilon,z})$, and give a rule to determine the highest
$\ell$-weight $M_{\varepsilon,z}$ which is a monomial in
$\mathbb{C}[\mathcal{Y}_{i,z}^{\pm1}]_{1\leq i\leq n,z\in\mathbb{C}^{\times}}$.

Exhibiting the $n$-dependence explicitly, let $P_n$ be the weight lattice of the simple Lie algebra of type $C_n$. Set
\begin{gather*}
E_n =\left\{\varepsilon=(\varepsilon_1,\ldots,\varepsilon_n)\mid\varepsilon_i=\pm1,1\leq i\leq n\right\}.
\end{gather*}
We def\/ine two weight functions
\begin{gather}
w_1 : \ E_n\longrightarrow P_n ,
\qquad \varepsilon\longmapsto\frac{1}{2}\sum^n\limits_{i=1}\varepsilon_i\epsilon_i ,
\label{weightfunction1}\\
w_2 : \ \{J\mid J\in\mathcal{J}^a (1\leq a\leq\sharp\mathcal{J}), \underset{\sim}{J}\subseteq\mathcal{J}\}\longrightarrow P_n ,
\qquad (j_1,\ldots,j_a)\longmapsto\sum^a\limits_{k=1}{\rm sgn}(j_k)\epsilon_{j_k} ,
\label{weightfunction2}
\end{gather}
where $\epsilon_0=0$, and
for $j\in\mathcal{J}$ we set
\begin{align*}
{\rm sgn}(j)=\begin{cases}
+,&j\preceq n,\\
0,&j=0,\\
-,&j\succeq\bar{n}.
\end{cases}
\end{align*}

Borrowing an idea from \cite{KS}, we introduce
monomials  $M_{\varepsilon,z}$ in
$\mathbb{C}[\mathcal{Y}_{i,z}^{\pm1}]_{1\leq i\leq n,z\in\mathbb{C}^{\times}}$
inductively as follows.
Def\/ine two operators  $\tau^{\mathcal{Y}}$, $\tau^z_{c}$ by
\begin{gather*}
\tau^{\mathcal{Y}}: \ \mathcal{Y}_{i,z}\rightarrow \mathcal{Y}_{i+1,z} ,\qquad
\tau^z_c: \ \mathcal{Y}_{i,z}\rightarrow \mathcal{Y}_{i,q^cz} .
\end{gather*}
For $n=1$ we set $M_{+,z}=\mathcal{Y}_{1,z}$ and $M_{-,z}=\mathcal{Y}^{-1}_{1,q^2z}$. In the general case we set
\begin{gather}
M_{(++\xi),z} =\tau^{\mathcal{Y}}(M_{(+\xi),z}) ,
\qquad M_{(+-\xi),z}=\mathcal{Y}_{1,q^{\frac{n}{2}}z}\tau^{\mathcal{Y}}(M_{(-\xi),z}) ,
\label{Mez1}
\\
M_{(-+\xi),z} =\mathcal{Y}^{-1}_{1,q^{\frac{n+2}{2}}z}\tau^z_1\tau^{\mathcal{Y}}(M_{(+\xi),z}) ,
\qquad M_{(--\xi),z}=\tau^z_1\tau^{\mathcal{Y}}(M_{(-\xi),z}) ,
\label{Mez2}
\end{gather}
where $\xi\in E_{n-2}$.

Now let $\varepsilon\in E_n$ with
$w_1(\varepsilon)=\frac{1}{2}\left(\sum^s_{k=1}\epsilon_{\varphi_k}-\sum^t_{k=1}\epsilon_{\psi_k}\right)$,
and consider the simple module $L(M_{\varepsilon,z})$.
We def\/ine $\mathcal{R}^{(n)}_{\varepsilon,z}$ and $\bar{\mathcal{R}}^{(n)}_{\varepsilon}$ through the
$q$-character and the character as follows.
\begin{gather*}
\mathcal{R}^{(n)}_{\varepsilon,z}=
 \chi_q\left(L(M_{\varepsilon,z})\right)\prod_{1\leq k\leq l\leq t}(1-x_{\psi_k}x_{\psi_l})
\prod_{1\leq k\leq s,1\leq l\leq t,\atop \varphi_k>\psi_l}\big(1-x_{\varphi_k}^{-1}x_{\psi_l}\big) ,\\
\bar{\mathcal{R}}^{(n)}_{\varepsilon}=
 \chi \left(L(M_{\varepsilon,z})\right)\prod_{1\leq k\leq l\leq t}(1-x_{\psi_k}x_{\psi_l})
\prod_{1\leq k\leq s,1\leq l\leq t,\atop \varphi_k>\psi_l}\big(1-x_{\varphi_k}^{-1}x_{\psi_l}\big) .
\end{gather*}
For the latter we have the following guess.

\begin{Conjecture}\label{spin-character-C_n}
For $\varepsilon=(\varepsilon_1,\ldots,\varepsilon_n)\in E_n$, we have
\begin{gather*}
\chi\left(L(M_{\varepsilon,z})\right)=
\prod_{1\leq i\leq n}\frac{1}{1-x_i^{-2\varepsilon_i}}
\prod_{1\leq i<j\leq n}\frac{1}{1-x_i^{-\varepsilon_i}x_j^{-\varepsilon_j}}.
\end{gather*}
\end{Conjecture}

In the special case $n=2$,
one can obtain the series for $C^{(1)}_2$ from that of $B^{(1)}_2$.
Let $\mathcal{Q}^{(1)}_{i,z}(1\leftrightarrow2)$ stand for the series for $B^{(1)}_2$ with $a=1$ given in Example~\ref{B12},
wherein we interchange $\mathcal{Y}_{1,z}$ with~$\mathcal{Y}_{2,z}$ and $A_{1,z}$ with $A_{2,z}$.
Then we have
\begin{gather}
\mathcal{R}^{(2)}_{(++),z}=
 \mathcal{Q}^{(1)}_{1,z}(1\leftrightarrow2) , \qquad
\mathcal{R}^{(2)}_{(+-),z}=\mathcal{Q}^{(1)}_{2,z}(1\leftrightarrow2) ,
\label{B2C2-1}\\
\mathcal{R}^{(2)}_{(-+),z}= \mathcal{Q}^{(1)}_{\bar{2},z}(1\leftrightarrow2) ,
\qquad \mathcal{R}^{(2)}_{(--),z}=\mathcal{Q}^{(1)}_{\bar{1},z}(1\leftrightarrow2) .
\label{B2C2-2}
\end{gather}

\subsection[Series $\mathcal{R}^{(n-1)}_{\varepsilon,z}$, $\mathcal{R}^{(n)}_{\varepsilon,z}$ for the spin nodes: case $D^{(1)}_n$]{Series $\boldsymbol{\mathcal{R}^{(n-1)}_{\varepsilon,z}$, $\mathcal{R}^{(n)}_{\varepsilon,z}}$ for the spin nodes: case $\boldsymbol{D^{(1)}_n}$}

For the two spin nodes of $D^{(1)}_n$, we follow the procedure for $C^{(1)}_n$.
We use the same weight functions \eqref{weightfunction1}, \eqref{weightfunction2} as for $C^{(1)}_n$.
We also introduce two subsets of $E_n$,
\begin{gather*}
E_{n,\varsigma} =\left\{\varepsilon\in E_n\mid\sharp\{\varepsilon_i=-1\}=\varsigma\pmod2\right\} ,\qquad  \varsigma=0,1 .
\end{gather*}
We def\/ine monomials $M_{\varepsilon,z}$ inductively by
\begin{gather}
M_{(++\xi),z}
=\tau^{\mathcal{Y}}(M_{(+\xi),z}) ,
\qquad M_{(+-\xi),z}=\mathcal{Y}_{1,q^{n-2}z}\tau^{\mathcal{Y}}(M_{(-\xi),z}) ,
\label{Mez3}
\\
M_{(-+\xi),z} =\mathcal{Y}^{-1}_{n,q^{l}z}\tau^z_2\tau^{\mathcal{Y}}(M_{(+\xi),z}) ,\qquad
M_{(--\xi),z}=\tau^z_2\tau^{\mathcal{Y}}(M_{(-\xi),z}) ,\label{Mez4}
\end{gather}
with $\xi\in E_{n-2}$ and the initial values
\begin{gather*}
M_{(++),z}=\mathcal{Y}_{2,z} ,\qquad M_{(--),z}=\mathcal{Y}^{-1}_{2,q^2z} ,\\
M_{(+-),z}=\mathcal{Y}_{1,z} ,\qquad M_{(-+),z}=\mathcal{Y}^{-1}_{1,q^2z} .
\end{gather*}

For $\varepsilon\in E_{n,\varsigma}$ ($\varsigma=0,1$)
with $w_1(\varepsilon)=\frac{1}{2}\left(\sum^s_{k=1}\epsilon_{\varphi_k}-\sum^t_{k=1}\epsilon_{\psi_k}\right)$, we def\/ine
\begin{gather*}
\mathcal{R}^{(n-\varsigma)}_{\varepsilon,z}= \chi_q\left(L(M_{\varepsilon,z})\right)\prod_{1\leq k< l\leq t}(1-x_{\psi_k}x_{\psi_l})
\prod_{1\leq k\leq s,1\leq l\leq t,\atop \varphi_k>\psi_l}\big(1-x_{\varphi_k}^{-1}x_{\psi_l}\big) ,\\
\bar{\mathcal{R}}^{(n-\varsigma)}_{\varepsilon}= \chi\left(L(M_{\varepsilon,z})\right)\prod_{1\leq k< l\leq t}(1-x_{\psi_k}x_{\psi_l})
\prod_{1\leq k\leq s,1\leq l\leq t,\atop \varphi_k>\psi_l}\big(1-x_{\varphi_k}^{-1}x_{\psi_l}\big) .
\end{gather*}

\begin{Conjecture}\label{spin-character-D_n}
For $\varepsilon=(\varepsilon_1,\ldots,\varepsilon_n)\in E_{n,\varsigma}$ $(\varsigma=0,1)$, we have
\begin{gather*}
\chi\left(L(M_{\varepsilon,z})\right)=\prod_{1\leq i<j\leq n}\frac{1}{1-x_i^{-\varepsilon_i}x_j^{-\varepsilon_j}} .
\end{gather*}
\end{Conjecture}

Although we do not have formulas for the series $\mathcal{R}^{(n-\varsigma)}_{\varepsilon,z}$ ($\varsigma=0,1$) for general $n$, in the special case $D^{(1)}_4$, similarly to $\mathcal{Q}^{(1)}_{i,z}$, we def\/ine $\mathcal{R}^{(3)}_{\varepsilon,z}$ $($resp.\ $\mathcal{R}^{(4)}_{\varepsilon,z})$ as certain limits of the q-characters of the KR modules $W^{(3)}_{m,z}$ $($resp. $W^{(4)}_{m,z}$). More precisely, one can obtain $\mathcal{R}^{(3)}_{\varepsilon,z}$ $($resp. $\mathcal{R}^{(4)}_{\varepsilon,z})$ from $\mathcal{Q}^{(1)}_{i,z}$ by exchanging $\mathcal{Y}_{1,z}$ with $\mathcal{Y}_{3,z}$
$($resp. $\mathcal{Y}_{4,z})$ and $A_{1,z}$ with $A_{3,z}$ $($resp.\ $A_{4,z})$. For example, we have
\begin{gather*}
\mathcal{Q}^{(1)}_{\bar{2},z}= \mathcal{Y}^{-1}_{2,q^5z}\mathcal{Y}_{1,q^4z}\sum_{k\geq0}\prod^k_{j=1}A^{-1}_{1,q^{6-2j}z} ,\qquad
\mathcal{Q}^{(1)}_{\bar{1},z}=\mathcal{Y}^{-1}_{1,q^6z} ,\\
\mathcal{R}^{(3)}_{(--+-),z}= \mathcal{Y}^{-1}_{2,q^5z}\mathcal{Y}_{3,q^4z}\sum_{k\geq0}\prod^k_{j=1}A^{-1}_{3,q^{6-2j}z} ,\qquad
\mathcal{R}^{(3)}_{(---+),z}=\mathcal{Y}^{-1}_{3,q^6z} ,\\
\mathcal{R}^{(4)}_{(--++),z}= \mathcal{Y}^{-1}_{2,q^5z}\mathcal{Y}_{4,q^4z}\sum_{k\geq0}\prod^k_{j=1}A^{-1}_{4,q^{6-2j}z} ,\qquad
\mathcal{R}^{(4)}_{(----),z}=\mathcal{Y}^{-1}_{4,q^6z} ,
\end{gather*}
and so on.

\section{Polynomial relations}\label{sec:polyrela}
In this section, we shall give the main results of this paper. We propose polynomial relations among the series def\/ined in Section~\ref{sec:series}, which are
expected from the $\psi$-system obtained in Section~\ref{sec:psi-system}.

We prepare some notation which we will use below.
For an element $j\in\mathcal{J}$, an element $J=(j_1,\ldots,j_a)\in\mathcal{J}^a$ and a positive integer $k$ with $1\leq k\leq a$, we set
\begin{gather*}
 \lambda_{J}=\sum^a_{l=1}\lambda_{j_l},\qquad
x_J=\prod^a_{l=1}x_{j_l},
\qquad
 J-j_k=(j_1,\ldots,\hat{j}_k,\ldots,j_a), \qquad (j,J)=(j,j_1,\ldots,j_a).
\end{gather*}
Let $\{v_i\}_{i\in\widetilde{\mathcal{J}}}$ be the standard basis of the vector representation
$V^{(1)}$ of ${}^L\mathfrak{g}$ (cf.\ Appendix~\ref{app1}).
Here $\widetilde{\mathcal{J}}=\{1,\ldots,n,0,\bar{0},\bar{n},\ldots,\bar{1}\}$ for $C^{(1)}_n$ and $\widetilde{\mathcal{J}}=\mathcal{J}$ for the other cases.
For $J=(j_1,\ldots,j_a)\in\widetilde{\mathcal{J}}^a$, we set $v_J=v_{j_1}\wedge v_{j_2}\wedge\cdots\wedge v_{j_a}$.

\subsection{Pl{\"u}cker-type relations}\label{subsection5.1}

We begin by writing down the $\psi$-system and the corresponding relations for
 $\mathcal{Q}^{(a)}_{J,z}$ for all nodes excepting the last (or the second last for $D^{(1)}_n$).
In fact they are just the consequence of the fact that
the $\mathcal{Q}^{(a)}_{J,z}$ are def\/ined by determinants of the  $\mathcal{Q}^{(1)}_{i,z}$'s.
So they are rather def\/initions, and we write them down just for uniformity reasons.
As we shall see later, the only non-trivial relations come from the last (or the second last) node.

For comparison, we start with the $\psi$-system.
Let $J_1=(i_1,\ldots,i_a)$, $J_2=(j_1,\ldots,j_a)$ be two elements of $\widetilde{\mathcal{J}}^a$.
The embedding $\iota$ \eqref{iotaA} of ${}^L\mathfrak{g}$-modules is explicitly given by
\begin{gather*}
\iota : \ \bigwedge^2 V^{(a)}_{\frac{1}{2t}} \hookrightarrow V^{(a-1)}\otimes V^{(a+1)},\\
\hphantom{\iota :}{} \ \ v_{J_1}\wedge v_{J_2} \mapsto\frac{1}{2}\left(\sum^a_{k=1}(-1)^{k-1}v_{J_1-i_k}\otimes v_{(i_k,J_2)}-\sum^a_{k=1}(-1)^{k-1}v_{J_2-j_k}\otimes v_{(j_k,J_1)}\right).
\end{gather*}
Therefore \eqref{psiA-1}, \eqref{psiB-1}, \eqref{psiC-1} and \eqref{psiD-1} imply
\begin{gather}
 Q^{(a)}_{J_1,-\frac{1}{2t}}(E)Q^{(a)}_{J_2,\frac{1}{2t}}(E)
-Q^{(a)}_{J_1,\frac{1}{2t}}(E)Q^{(a)}_{J_2,-\frac{1}{2t}}(E)\omega^{\lambda_{J_1}-\lambda_{J_2}}
\nonumber\\
\qquad{} = \sum^{a}_{k=1}(-1)^k
Q^{(a-1)}_{J_2-j_k}(E)Q^{(a+1)}_{(j_k,J_1)}(E)\omega^{-\lambda_{j_k}}.\label{TypeA-relation}
\end{gather}

The corresponding relations for $\mathcal{Q}^{(a)}_{J,z}$  read as follows.

\begin{Proposition}
For two increasing elements $J_1=(i_1,i_2,\ldots,i_{a})$ and $J_2=(j_1,j_2,\ldots,j_{a})$ of~$\mathcal{J}^a$, we have
\begin{gather}
 \mathcal{Q}^{(a)}_{J_1,q_1^{-1}z}\mathcal{Q}^{(a)}_{J_2,q_1z}
-\mathcal{Q}^{(a)}_{J_1,q_1z}\mathcal{Q}^{(a)}_{J_2,q_1^{-1}z}\frac{x_{J_1}}{x_{J_2}}
=\sum^{a}_{k=1}(-1)^k
\mathcal{Q}^{(a-1)}_{J_2-j_k,z}\mathcal{Q}^{(a+1)}_{(j_k,J_1),z}x_{j_k}^{-1}.\label{Plucker}
\end{gather}
\end{Proposition}

This is an immediate consequence of the Sylvester identity
(see e.g.~\cite[p.~108]{Fulton}). Note that~\eqref{TypeA-relation} and~\eqref{Plucker}
has the same form under the identif\/ication
\begin{gather*}
Q^{(a)}_{J}(E)\leftrightarrow\mathcal{Q}^{(a)}_{J,z} ,\qquad E\leftrightarrow z ,
\qquad\Omega^{\frac{1}{2t}}\leftrightarrow q_1 ,\qquad\omega^{\lambda_i}\leftrightarrow x_i .
\end{gather*}

\subsection[Wronskian identity for $A_n^{(1)}$]{Wronskian identity for $\boldsymbol{A_n^{(1)}}$}
For type $A_n^{(1)}$, the non-trivial relation is the following
`Wronskian identity'.

\begin{Theorem}\label{QQ-A_n}
Let $\mathfrak{g}=A^{(1)}_n$. We have
\begin{gather*}
\det\left(\mathcal{Q}^{(1)}_{\nu,q^{n+2-2\mu}z}x^{-\mu+\nu}_{\nu}\right)^{n+1}_{\mu,\nu=1}=1.
\end{gather*}
\end{Theorem}

Note that the left-hand side is just
$\mathcal{Q}^{(n+1)}_{J,z}(-1)^{(n+1)n/2}\!\times\!\prod^{n}_{b=1}x_b^{n+1-b}$ with $J=(1,2,\ldots,n{+}1)$.
We prove Theorem~\ref{QQ-A_n} in Appendix~\ref{app2}.

\subsection[Polynomial relations related to the last node: case $B^{(1)}_n$]{Polynomial relations related to the last node: case $\boldsymbol{B^{(1)}_n}$}\label{subsection5.3}

In this subsection we
let $\mathfrak{g}=B^{(1)}_n$.
We give conjectural relations which correspond to the last identity~\eqref{psiB-2} of the $\psi$-system.
We f\/irst give identities for the connection coef\/f\/icients $Q^{(a)}_{J}(E)$.

In this subsection, for an element $j\in\mathcal{J}$ and $J=(j_1,\ldots,j_a)\in\mathcal{J}^a$ we def\/ine
\begin{gather*}
c(j)=\begin{cases}
(-1)^j,&j\preceq n,\\
(-1)^{\bar{j}-1},&j\succeq\bar{n},
\end{cases} \qquad c(J)=\prod^a_{k=1}c(j_k).
\end{gather*}
We also denote by $\sigma(J,J^*)$ the signature of the permutation $(1,\ldots,n,\bar{n},\ldots,\bar{1})\mapsto(J,J^*)$.

The bilinear form $\langle~,~\rangle: V\times V_{\frac{1}{2}}\to \mathbb{C}$ given by
$\langle v_i,v_j\rangle=c(i)\delta_{i,\bar j}$
is ${}^L\mathfrak{g}$-invariant.
From this it is easy to see that we have the following ${}^L\mathfrak{g}$-module isomorphism
\begin{gather*}
V^{(a)}_{\frac{1}{2}} \simeq V^{(2n-a)}, \qquad  a=1,\ldots,n, \qquad v_J \mapsto c(J)\sigma(J,J^*)v_{\bar{J}^*},
\end{gather*}
which implies
\begin{gather}
Q^{(a)}_{J,\frac{1}{2}}(E)=c(J)\sigma(J,J^*)Q^{(2n-a)}_{\overline{J^*}}(E).\label{QQB}
\end{gather}
Let $J_1=(i_1,\ldots,i_n)$, $J_2=(j_1,\ldots,j_n)$ be two elements of $\mathcal{J}^n$. Then the relation \eqref{TypeA-relation} together with  \eqref{QQB} imply
\begin{gather}
 Q^{(n)}_{\overline{J^*_1}}(E)Q^{(n)}_{J_2,\frac{1}{2}}(E)c\left(\overline{J^*_1}\right)\sigma\left(\overline{J^*_1},\overline{J_1}\right)
\omega^{-\lambda_{J_1}}-Q^{(n)}_{\overline{J^*_2}}(E)Q^{(n)}_{J_1,\frac{1}{2}}(E)
c\left(\overline{J^*_2}\right)\sigma\left(\overline{J^*_2},\overline{J_2}\right)
\omega^{-\lambda_{J_2}}
\nonumber\\
\qquad{}= \sum^n_{k=1}(-1)^kQ^{(n-1)}_{J_2-j_k}(E)Q^{(n-1)}_{\overline{J^*_1}-\overline{j_k},\frac{1}{2}}(E)
c\left(\overline{(j_k,J_1)}^*\right)\sigma\left(\overline{(j_k,J_1)}^*,\overline{(j_k,J_1)}\right)\omega^{-\lambda_{J_1}-\lambda_{j_k}}.\!\!\!\label{conn-coef-B_n}
\end{gather}

For a given element $i\in\mathcal{J}$ and $J\in\mathcal{J}^a$, let
\begin{gather}
X_J=\prod_{\mu,\nu\in J,\atop \mu\prec \nu}\big({-}x_{\mu}^{\frac{1}{2}\delta_{\mu,\bar{\nu}}-1}\big), \qquad x_{i,J}=\prod_{j\in J,\atop i\prec j}(-x_i^{1-\frac{1}{2}\delta_{j,\bar{i}}})\prod_{j\in J^*,\atop j\prec i}\big(x_j^{\frac{1}{2}\delta_{i,\bar{j}}-1}\big).\label{notation-B_n}
\end{gather}
The counterpart of \eqref{conn-coef-B_n} for the series $\mathcal{Q}^{(a)}_{J,z}$ is given by the following conjecture.

\begin{Conjecture} \label{conj-B_n}
Let $J_1=(i_1,\ldots,i_n), J_2=(j_1,\ldots,j_n)\in \mathcal{J}^n$ be increasing.
Then we have
\begin{gather*}
 \mathcal{Q}^{(n)}_{\overline{J^*_1},z}\mathcal{Q}^{(n)}_{J_2,qz}X_{\overline{J^*_1}}X_{J_2}x^{-\frac{1}{2}}_{J_1}
-\mathcal{Q}^{(n)}_{\overline{J^*_2},z}\mathcal{Q}^{(n)}_{J_1,qz}X_{\overline{J^*_2}}X_{J_1}
x_{J_2}^{-\frac{1}{2}}\\
\qquad{}= -\sum_{j_k\in J_2\cap J_1^*, \atop 1\leq k\leq n}\mathcal{Q}^{(n-1)}_{J_2-j_k,z}\mathcal{Q}^{(n-1)}_{\overline{J^*_1}-\overline{j_k},qz}X_{J_2-j_k}X_{\overline{J^*_1}-\overline{j_k}}
\frac{x_{j_k,J_1}}{x_{j_k,J_2}}x^{-\frac{1}{2}}_{J_1}x^{-\frac{1}{2}}_{j_k}.
\end{gather*}
\end{Conjecture}

When $n=2$,  we have verif\/ied the above conjecture  by direct computations.
To save space we do not write the proofs here.
We have also checked them for $n=3$ by Mathematica~5.0 up to certain degree, counting the degree of $A_{a,q^kz}$ to be $k$.  For $n\geq 4$, although at the moment we do not have a proof of these identities, in Appendix \ref{app3}, we show that the identities hold when specialized to the characters.

\begin{Remark}
It is worth noting that, although the connection coef\/f\/icients satisfy the linear relation \eqref{QQB},
there are no analogs for the series $\mathcal{Q}^{(a)}_{J,z}$.
For example, in the simplest case $B_2^{(1)}$,
a direct computation gives
$\mathcal{Q}^{(3)}_{(2\bar{2}\bar{1}),z}=-x_2\mathcal{Y}^{-1}_{1,q^4z}\mathcal{Y}^{-1}_{2,q^2z}\mathcal{Y}_{2,q^3z}
-\mathcal{Y}^{-1}_{1,q^2z}\mathcal{Y}_{2,q^3z}\mathcal{Y}^{-1}_{2,q^4z}$, while $\mathcal{Q}^{(1)}_{\bar{1},z}
=\mathcal{Y}^{-1}_{1,q^3z}$, and
there is no function $f(x_1,x_2)$ of $x_1$, $x_2$
such that $\mathcal{Q}^{(3)}_{(2\bar{2}\bar{1}),z}=f(x_1,x_2)\mathcal{Q}^{(1)}_{\bar{1},z}$.
So only the quadratic relations survive.
\end{Remark}

\subsection[Polynomial relations related to spin node: case $C^{(1)}_n$]{Polynomial relations related to spin node: case $\boldsymbol{C^{(1)}_n}$}

In this subsection, $\mathfrak{g}=C^{(1)}_n$.
We shall give conjectural relations for the series $\mathcal{Q}^{(a)}_{J,z}$ and $\mathcal{R}^{(n)}_{\varepsilon,z}$.

For an increasing element $J\in\mathcal{J}^a$
we can write
\begin{gather*}
 \underset{\sim}{J}=\left\{\alpha_1,\ldots,\alpha_m,\bar{\beta}_1,\ldots,\bar{\beta}_r,\sigma_1,\ldots,\sigma_u, \bar{\sigma}_1,\ldots,\bar{\sigma}_u\right\}\cup(\underset{\sim}{J}\cap\{0\}) ,\\
 \underset{\sim}{J^*}=\left\{\bar{\alpha}_1,\ldots,\bar{\alpha}_m,\beta_1,\ldots,\beta_r,\eta_1,\ldots,\eta_{v}, \bar{\eta}_1,\ldots,\bar{\eta}_{v}\right\}\cup(\underset{\sim}{J^*}\cap\{0\}) ,
\end{gather*}
where $\alpha_i$, $\beta_i$, $\sigma_i$, $\eta_i$ are mutually distinct elements of
$\{1,\ldots,n\}$.
Let $\pmb{\sigma}=\{\sigma_1,\ldots,\sigma_u\}$, $\pmb{\eta}=\{\eta_1,\ldots,\eta_v\}$. Note that
\begin{gather*}
v=\begin{cases}
u+2,&0\in\underset{\sim}{J}\;\text{and}\;J\in\mathcal{J}^{n-1},\\
u+1,&0\notin\underset{\sim}{J}\;\text{and}\;J\in\mathcal{J}^{n-1},\;\text{or}\;0\in\underset{\sim}{J}\;\text{and}\;J\in\mathcal{J}^{n},\\
u,&0\notin\underset{\sim}{J}\;\text{and}\;J\in\mathcal{J}^{n}.
\end{cases}
\end{gather*}

Let $J\in\mathcal{J}^a$ be increasing. We introduce the following condition for $\varepsilon, \varepsilon'\in E_n$:

{\noindent\bf\underline{Condition $C_J$}:}
\begin{itemize}\itemsep=0pt
  \item $w_1(\varepsilon)+w_1(\varepsilon')=w_2(J)$,
  \item if $w_1(\varepsilon)-w_1(\varepsilon')=\sum^s\limits_{k=1}\epsilon_{\gamma_k}-\sum^t\limits_{k=1}\epsilon_{\delta_k}$, then $\pmb{\sigma}\subseteq\pmb{\gamma}$.
\end{itemize}
Here $w_1$, $w_2$ are the weight functions def\/ined by \eqref{weightfunction1}, \eqref{weightfunction2} and $\pmb{\gamma}=\{\gamma_1,\ldots,\gamma_s\}$, $\pmb{\delta}=\{\delta_1,\ldots,\delta_t\}$. It is easy to see $\pmb{\sigma}\cup\pmb{\eta}=\pmb{\gamma}\cup\pmb{\delta}$.

In the following Conjectures \ref{conj-C_n-1}--\ref{conj-C_n-4},
$\mathfrak{g}=C^{(1)}_n$, $\pmb{\sigma}$, $\pmb{\eta}$ are disjoint subsets of
$\{1,2,\ldots,n\}$.

\begin{Conjecture}\label{conj-C_n-1}
Suppose $J\in\mathcal{J}^{n-1}$ is increasing and $0\notin\underset{\sim}{J}$. Then we have
\begin{gather*}
 (-1)^{\frac{u(u-1)}{2}}\mathcal{Q}^{(n-1)}_{J,z}x_{\pmb{\sigma}}^{-2u}x_{\pmb{\eta}}^{2-u}\prod_{i,j\in J,\atop i\prec j}\big({-}x_i^{-1}\big)\prod_{1\leq k\leq u,1\leq l\leq u+1,\atop \sigma_k\succ\eta_l}\left(-\frac{x_{\sigma_k}}{x_{\eta_l}}\right)\\
\qquad{}= \sum_{}(-1)^{\frac{t(t+1)}{2}}\mathcal{R}^{(n)}_{\varepsilon,qz}\mathcal{R}^{(n)}_{\varepsilon',q^{-1}z}
x_{\pmb{\gamma}}^{2-2u}x_{\pmb{\delta}}^{1-t}\prod_{1\leq k\leq s,1\leq l\leq t,\atop \gamma_k\succ\delta_l}\left(-\frac{x_{\gamma_k}}{x_{\delta_l}}\right),
\end{gather*}
where the sum is taken over all $\varepsilon, \varepsilon'\in E_n$ satisfying \textup{Condition $C_J$}
where $\pmb{\gamma}$, $\pmb{\delta}$ are determined there.
\end{Conjecture}

\begin{Conjecture}\label{conj-C_n-2}
Suppose $J\in\mathcal{J}^{n-1}$ is increasing and $0\in\underset{\sim}{J}$. Then we have
\begin{gather*}
 (-1)^{\frac{u(u-1)}{2}}\mathcal{Q}^{(n-1)}_{J,z}x_{\pmb{\sigma}}^{-2u}x_{\pmb{\eta}}^{2-u}\prod_{i,j\in J\atop i\prec j}(-x_i^{-1})\prod_{1\leq k\leq u,1\leq l\leq u+2,\atop \sigma_k\succ\eta_l}\left(-\frac{x_{\sigma_k}}{x_{\eta_l}}\right)\\
\qquad{}= \sum_{}(-1)^{\frac{t(t-1)}{2}}\mathcal{R}^{(n)}_{\varepsilon,qz}\mathcal{R}^{(n)}_{\varepsilon',q^{-1}z}
x_{\pmb{\gamma}}^{2-2u}x_{\pmb{\delta}}^{2-t}\prod_{1\leq k\leq s,1\leq l\leq t,\atop \gamma_k\succ\delta_l}\left(-\frac{x_{\gamma_k}}{x_{\delta_l}}\right),
\end{gather*}
where the sum is taken over all $\varepsilon, \varepsilon'\in E_n$ satisfying \textup{Condition $C_J$}
where $\pmb{\gamma}$, $\pmb{\delta}$ are determined there.
\end{Conjecture}

\begin{Conjecture}\label{conj-C_n-3}
Suppose $J\in\mathcal{J}^{n}$ is increasing and $0\notin\underset{\sim}{J}$. Then we have
\begin{gather*}
 (-1)^{\frac{u(u-1)}{2}}\mathcal{Q}^{(n)}_{J,z}x_{\pmb{\sigma}}^{-2u}x_{\pmb{\eta}}^{1-u}\prod_{i,j\in J,\atop i\prec j}\big({-}x_i^{-1}\big)\prod_{1\leq k\leq u,1\leq l\leq u,\atop \sigma_k\succ\eta_l}\left(-\frac{x_{\sigma_k}}{x_{\eta_l}}\right)\\
\qquad{} = \sum_{}(-1)^{\frac{s(s-1)}{2}}\mathcal{R}^{(n)}_{\varepsilon,q^{\frac{1}{2}}z}\mathcal{R}^{(n)}_{\varepsilon',q^{-\frac{1}{2}}z}
x_{\pmb{\gamma}}^{1-2u}x_{\pmb{\delta}}^{-t}\prod_{1\leq k\leq s,1\leq l\leq t,\atop \gamma_k\succ\delta_l}\left(-\frac{x_{\gamma_k}}{x_{\delta_l}}\right),
\end{gather*}
where the sum is taken over all $\varepsilon, \varepsilon'\in E_n$ satisfying \textup{Condition $C_J$}
where $\pmb{\gamma}$, $\pmb{\delta}$ are determined there.
\end{Conjecture}

\begin{Conjecture}\label{conj-C_n-4}
Suppose $J\in\mathcal{J}^{n}$ is increasing and $0\in\underset{\sim}{J}$. Then we have
\begin{gather*}
 (-1)^{\frac{u(u-1)}{2}}\mathcal{Q}^{(n)}_{J,z}x_{\pmb{\sigma}}^{-2u}x_{\pmb{\eta}}^{1-u}\prod_{i,j\in J,\atop i\prec j}\big({-}x_i^{-1}\big)\prod_{1\leq k\leq u,1\leq l\leq u+1,\atop \sigma_k\succ\eta_l}\left(-\frac{x_{\sigma_k}}{x_{\eta_l}}\right)\\
\qquad{} = \sum_{}(-1)^{\frac{s(s-1)}{2}}\mathcal{R}^{(n)}_{\varepsilon,q^{\frac{1}{2}}z}\mathcal{R}^{(n)}_{\varepsilon',q^{-\frac{1}{2}}z}
x_{\pmb{\gamma}}^{1-2u}x_{\pmb{\delta}}^{1-t}\prod_{1\leq k\leq s,1\leq l\leq t,\atop \gamma_k\succ\delta_l}\left(-\frac{x_{\gamma_k}}{x_{\delta_l}}\right),
\end{gather*}
where the sum is taken over all $\varepsilon, \varepsilon'\in E_n$ satisfying \textup{Condition $C_J$}
where $\pmb{\gamma}$, $\pmb{\delta}$ are determined there.
\end{Conjecture}

For $n=2$ one can verify these relations using the relation between
 $C^{(1)}_2$ and $B^{(1)}_2$ (see \eqref{B2C2-1}, \eqref{B2C2-2}).
For general $n$, we show in Appendix \ref{app4} the validity of the identities when specialized to
the ordinary characters $\bar{\mathcal{Q}}^{(a)}_J$ and $\bar{\mathcal{R}}^{(n)}_{\varepsilon}$.

\begin{Remark}\label{Remark-C_n}
In the conjectured identities given above,
$\widehat{Q}^{(1)}_{0,z}$ and $\widehat{Q}^{(1)}_{\bar{0},z}$ enter only through the sum
$\mathcal{Q}^{(1)}_{0,z}=\widehat{Q}^{(1)}_{0,z}+\widehat{Q}^{(1)}_{\bar{0},z}$.
However, since the weight $0$ of $V^{(1)}$ has multiplicity $2$,
it is more natural
to consider $\widehat{Q}^{(1)}_{0,z}$ and $\widehat{Q}^{(1)}_{\bar{0},z}$ separately.
The $\psi$-system also suggests that there are
identities involving them
separately.
This is indeed the case for $n=2$
where we have, for example,
\begin{gather*}
 \mathcal{Q}^{(1)}_{1,q^{-\frac{1}{2}}z}\widehat{\mathcal{Q}}^{(1)}_{0,q^{\frac{1}{2}}z}
-x_1\widehat{\mathcal{Q}}^{(1)}_{\bar{0},q^{-\frac{1}{2}}z}\mathcal{Q}^{(1)}_{1,q^{\frac{1}{2}}z}
=-\frac{x_1}{x_2}\mathcal{R}^{(2)}_{(++),q^{-\frac{1}{2}}z}\mathcal{R}^{(2)}_{(+-),q^{\frac{1}{2}}z},\\
 \mathcal{Q}^{(1)}_{1,q^{-\frac{1}{2}}z}\widehat{\mathcal{Q}}^{(1)}_{\bar{0},q^{\frac{1}{2}}z}
-x_1\widehat{\mathcal{Q}}^{(1)}_{0,q^{-\frac{1}{2}}z}\mathcal{Q}^{(1)}_{1,q^{\frac{1}{2}}z}
=-x_1\mathcal{R}^{(2)}_{(+-),q^{-\frac{1}{2}}z}\mathcal{R}^{(2)}_{(++),q^{\frac{1}{2}}z}.
\end{gather*}
On the other hand,
if we def\/ine $\mathcal{Q}^{(a)}_{J,z}$ as determinants treating $0$ and $\bar0$ independently, then when~$J$ contains both~$0$ and $\bar 0$ we have  $\bar{\mathcal{Q}}^{(a)}_J=0$
because $x_0=x_{\bar0}=1$.
Hence $\mathcal{Q}^{(a)}_{J,z}$ can not be explained as a $q$-character,
so we must modify the working hypothesis.
At the moment we do not know how to do that.
\end{Remark}

\subsection[Polynomial relations related to spin nodes: case $D^{(1)}_n$]{Polynomial relations related to spin nodes: case $\boldsymbol{D^{(1)}_n}$}

Finally we consider the case $\mathfrak{g}=D^{(1)}_n$.
As in the case of  $C^{(1)}_n$,
for an increasing element $J\in\mathcal{J}^a$ we can write
\begin{gather*}
 \underset{\sim}{J}=\left\{\alpha_1,\ldots,\alpha_m,\bar{\beta}_1,\ldots,\bar{\beta}_r,\sigma_1,\ldots,\sigma_u, \bar{\sigma}_1,\ldots,\bar{\sigma}_u\right\},\\
 \underset{\sim}{J^*}=\left\{\bar{\alpha}_1,\ldots,\bar{\alpha}_m,\beta_1,\ldots,\beta_r,\eta_1,\ldots,\eta_{v}, \bar{\eta}_1,\ldots,\bar{\eta}_{v}\right\},
\end{gather*}
with mutually distinct $\alpha_i,\beta_i,\sigma_i,\eta_i\in \{1,\ldots,n\}$.
Let $\pmb{\sigma}=\{\sigma_1,\ldots,\sigma_u\}$, $\pmb{\eta}=\{\eta_1,\ldots,\eta_v\}$. Note that
\begin{gather*}
v=\begin{cases}
u+2,&J\in\mathcal{J}^{n-2},\\
u+1,&J\in\mathcal{J}^{n-1}.
\end{cases}
\end{gather*}

Let $J\in\mathcal{J}^a$ be increasing. We introduce the following conditions for $\varepsilon, \varepsilon'\in E_n$:

{\noindent\bf\underline{Condition $D_{J,\varsigma}$}:}
\begin{itemize}\itemsep=0pt
  \item $w_1(\varepsilon)+w_1(\varepsilon')=w_2(J)$,
  \item if $w_1(\varepsilon)-w_1(\varepsilon')
=\sum^s\limits_{k=1}\epsilon_{\gamma_k}-\sum^t\limits_{k=1}\epsilon_{\delta_k}$,
then $\pmb{\sigma}\subseteq\pmb{\gamma}$, and
$t\equiv r+\varsigma\pmod2$.
\end{itemize}
Here $\varsigma\in\{0,1\}$.

{\noindent\bf\underline{Condition $D_J$}:}
\begin{itemize}\itemsep=0pt
  \item $w_1(\varepsilon)+w_1(\varepsilon')=w_2(J)$,
  \item if $w_1(\varepsilon)-w_1(\varepsilon')=\sum^s\limits_{k=1}\epsilon_{\gamma_k}-\sum^t\limits_{k=1}\epsilon_{\delta_k}$, then $\pmb{\sigma}\subseteq\pmb{\gamma}$, and $t\equiv r\pmod2$.
\end{itemize}

Here $w_1$, $w_2$ are the weight functions def\/ined by \eqref{weightfunction1}, \eqref{weightfunction2}. We set $\pmb{\gamma}=\{\gamma_1,\ldots,\gamma_s\}$, $\pmb{\delta}=\{\delta_1,\ldots,\delta_t\}$. It is easy to see $\pmb{\sigma}\cup\pmb{\eta}=\pmb{\gamma}\cup\pmb{\delta}$.

In the following Conjectures \ref{conj-D_n-1}, \ref{conj-D_n-2},
$\mathfrak{g}=D^{(1)}_n$, and $\pmb{\sigma}$ and $\pmb{\eta}$ are disjoint subsets of $\{1,\ldots,n\}$.

\begin{Conjecture}\label{conj-D_n-1}
Let $J$ be an increasing element of $\mathcal{J}^{n-2}$, and let $\varsigma=0,1$.
Then we have
\begin{gather*}
 (-1)^{\frac{u(u-1)}{2}}\mathcal{Q}^{(n-2)}_{J,z}x_{\pmb{\sigma}}^{1-2u}x_{\pmb{\eta}}^{1-u}\prod_{i,j\neq\bar{i}\in J,\atop i\prec j}\big({-}x_i^{-1}\big)\prod_{1\leq k\leq u,1\leq l\leq u+2,\atop \sigma_k\succ\eta_l}\left(-\frac{x_{\sigma_k}}{x_{\eta_l}}\right)\\
\qquad{} = \sum_{}(-1)^{\frac{s(s+1)}{2}}\mathcal{R}^{(n-\varsigma)}_{\varepsilon,qz}\mathcal{R}^{(n-\varsigma)}_{\varepsilon',q^{-1}z}
x_{\pmb{\gamma}}^{-2u}x_{\pmb{\delta}}^{3-t}\prod_{1\leq k\leq s,1\leq l\leq t,\atop \gamma_k\succ\delta_l}\left(-\frac{x_{\gamma_k}}{x_{\delta_l}}\right),
\end{gather*}
where the sum is taken over all $\varepsilon, \varepsilon'\in E_{n,\varsigma}$ satisfying \textup{Condition $D_{J,\varsigma}$}, where $\pmb{\gamma}$ and $\pmb{\delta}$ are determined there.
\end{Conjecture}

\begin{Conjecture}\label{conj-D_n-2}
Let $J$ be an increasing element of $\mathcal{J}^{n-1}$. Then we have
\begin{gather*}
 (-1)^{\frac{u(u-1)}{2}}\mathcal{Q}^{(n-1)}_{J,z}x_{\pmb{\sigma}}^{-1-2u}x_{\pmb{\eta}}^{-u}\prod_{i,j\neq\bar{i}\in J,\atop i\prec j}\big({-}x_i^{-1}\big)\prod_{1\leq k\leq u,1\leq l\leq u+1,\atop \sigma_k\succ\eta_l}\left(-\frac{x_{\sigma_k}}{x_{\eta_l}}\right)\\
\qquad{}= \sum_{}(-1)^{\frac{t(t-1)}{2}}\mathcal{R}^{(n-1)}_{\varepsilon,z}\mathcal{R}^{(n)}_{\varepsilon',z}
x_{\pmb{\gamma}}^{-2u}x_{\pmb{\delta}}^{1-t}\prod_{1\leq k\leq s,1\leq l\leq t,\atop \gamma_k\succ\delta_l}\left(-\frac{x_{\gamma_k}}{x_{\delta_l}}\right),
\end{gather*}
where the sum is taken over all $\varepsilon\in E_{n,1}$, $\varepsilon'\in E_{n,0}$ satisfying \textup{Condition $D_J$}, where $\pmb{\gamma}$ and $\pmb{\delta}$ are determined there.
\end{Conjecture}

For $n=4$, we have checked by Mathematica 5.0 the conjectures hold up to certain degree.
For general $n$, we prove  in Appendix \ref{app4} the conjectures specialized to the ordinary
characters~$\bar{\mathcal{Q}}^{(a)}_J$ and $\bar{\mathcal{R}}^{(n-1)}_{\varepsilon}$,
$\bar{\mathcal{R}}^{(n)}_{\varepsilon}$.

\section{Conclusion}\label{sec:conclusion}

In this section, we give a summary of this paper.
We have done the following things:
\begin{enumerate}\itemsep=0pt
  \item To each weight of the fundamental representation $V^{(a)}$, we associated a formal series $\mathcal{Q}^{(a)}_{J,z}$ or $\mathcal{R}^{(a)}_{\varepsilon,z}$. We expect that with some simple factors the formal series are $q$-characters of certain irreducible modules of $U_q(\mathfrak{b})$.
  \item Under suitable identif\/ications, using relations for the connection coef\/f\/icients implied by the $\psi$-system, we proposed the following conjecture relations for the series $\mathcal{Q}^{(a)}_{J,z}$ and $\mathcal{R}^{(a)}_{\varepsilon,z}$:
      \begin{itemize}\itemsep=0pt
        \item $\mathfrak{g}=A^{(1)}_n, B^{(1)}_n, C^{(1)}_n, D^{(1)}_n$, Proposition~\ref{TypeA-relation}. This is the Pl\"{u}cker type relations which is in fact the def\/ining relations of the series $\mathcal{Q}^{(a)}_{J,z}$;
        \item $\mathfrak{g}=A^{(1)}_n$, Theorem \ref{Plucker}. This is the Wronskian identity for the $q$-characters of $U_q(\mathfrak{b})$ to be proved in Appendix \ref{app2}.
        \item $\mathfrak{g}=B^{(1)}_n$, Conjecture \ref{conj-B_n}.
        \item $\mathfrak{g}=C^{(1)}_n$, Conjectures \ref{conj-C_n-1}--\ref{conj-C_n-4}.
        \item $\mathfrak{g}=D^{(1)}_n$, Conjectures \ref{conj-D_n-1}, \ref{conj-D_n-2}.
      \end{itemize}
  For the last three cases, we support our conjectures by checking the following.
  For the special cases $\mathfrak{g}=B^{(1)}_2$, $\mathfrak{g}=C^{(1)}_2$,  we proved the conjectures by direct computations. For the  cases $\mathfrak{g}=B^{(1)}_3$ and $\mathfrak{g}=D^{(1)}_4$, we checked the conjectures up to some degrees by Mathematica. When specialized to characters, the conjectured relations hold in all cases.
This will be proved in Appendix \ref{app3} or \ref{app4}.
\end{enumerate}
We hope we have presented reasonable grounds to suggest that the correspondence between
the connection coef\/f\/icients of certain dif\/ferential equations
and the $q$-characters of the
Borel subalgebra $U_q(\mathfrak{b})$
supplies an ef\/fective way to f\/ind polynomial relations. Certainly this is only the f\/irst step, and
more serious checks are desirable along with attempts toward proving these identities.

\appendix

\section[Vector representation of ${}^L\mathfrak{g}$]{Vector representation of $\boldsymbol{{}^L\mathfrak{g}}$}\label{app1}

Following \cite[Appendix 2]{DS},
we give an explicit realization of the vector representation of ${}^L\mathfrak{g}$.
The symbol $E_{i,j}$ stands for the matrix unit
$\bigl(\delta_{i,a}\delta_{j,b}\bigr)_{a,b=1,\ldots,N}$ where $N=\dim V^{(1)}$.

\medskip

\noindent{\underline{${}^L\mathfrak{g}=A^{(1)}_n$}}:\quad
\begin{gather*}
 e_0=E_{n+1,1} ,\quad e_i=E_{i,i+1},\quad  1\le i\le n ,
\\
 f_0=E_{1,n+1} ,\quad f_i=E_{i+1,i},\quad  1\le i\le n ,
\\
 h_0=-E_{1,1}+E_{n+1,n+1} ,\quad h_i=E_{i,i}-E_{i+1,i+1},\quad  1\le i\le n ,
\\
 h_\rho=\diag\left(-\frac{n}{2},\dots,\frac{n}{2}\right) ,
\quad
 \ell=\diag\bigl(\mu_1,\dots,\mu_{n+1}\bigr), \quad  \sum_{i=1}^{n+1}\mu_i=0 .
\end{gather*}

\noindent{\underline{${}^L\mathfrak{g}=A^{(2)}_{2n-1}$}}:
\begin{gather*}
 e_0=\frac{1}{2}(E_{1,2n-1}+E_{2,2n}) ,\quad
e_i=E_{i+1,i}+E_{2n+1-i,2n-i},\quad  1\le i\le n-1 ,\quad
 e_n=E_{n+1,n} ,
\\
 f_0=2(E_{2n-1,1}+E_{2n,2}),\quad
f_i=E_{i,i+1}+E_{2n-i,2n+1-i},\quad  1\le i\le n-1,\quad
 f_n=E_{n,n+1},
\\
 h_0=E_{1,1}+E_{2,2}-E_{2n-1,2n-1}-E_{2n,2n},\\
 h_i=-E_{i,i}+E_{i+1,i+1}-E_{2n-i,2n-i}+E_{2n+1-i,2n+1-i},\quad  1\le i\le n-1,\\
 h_n=-E_{n,n}+E_{n+1,n+1},\\
 h_\rho=\diag\left(-n+\frac{1}{2},\dots,-\frac{1}{2},\frac{1}{2},\dots,n-\frac{1}{2}\right) ,\quad
 \ell=\diag\bigl(\mu_1,\dots,\mu_n,-\mu_n,\dots,-\mu_1\bigr) .
\end{gather*}

\noindent{\underline{${}^L\mathfrak{g}=D^{(2)}_{n+1}$}}:
\begin{gather*}
 e_0=E_{n+2,1}+E_{2n+2,n+2},\quad
e_i=E_{i,i+1}+E_{2n+2-i,2n+3-i},\quad  1\le i\le n-1 ,\\
 e_n=2(E_{n,n+1}+E_{n+1,n+3}),\\
 f_0=2(E_{1,n+2}+E_{n+2,2n+2}),\quad
f_i=E_{i+1,i}+E_{2n+3-i,2n+2-i},\quad 1\le i\le n-1 ,\\
 f_n=E_{n+1,n}+E_{n+3,n+1},\\
 h_0=2(-E_{1,1}+E_{2n+2,2n+2}),\\
 h_i=E_{i,i}-E_{i+1,i+1}+E_{2n+2-i,2n+2-i}-E_{2n+3-i,2n+3-i},\quad 1\le i\le n-1,\\
 h_n=2(E_{n,n}-E_{n+3,n+3}),\\
 h_\rho=\diag\bigl({-}n,\dots,-1,0,0,1,\dots,n\bigr),
\quad
\ell=\diag\bigl(\mu_1,\dots,\mu_n,0,0,-\mu_n,\dots,-\mu_1\bigr).
\end{gather*}

\noindent{\underline{${}^L\mathfrak{g}=D^{(1)}_{n}$}}:
\begin{gather*}
 e_0=\frac{1}{2}(E_{2n-1,1}+E_{2n,2}),\quad
e_i=E_{i,i+1}+E_{2n-i,2n+1-i},\quad  1\le i\le n-1,\\
e_n=2(E_{n-1,n+1}+E_{n,n+2}),
\\
 f_0=2\bigl(E_{1,2n-1}+E_{2,2n}\bigr),\quad
f_i=E_{i+1,i}+E_{2n+1-i,2n-i},\quad 1\le i\le n-1 ,\\
 f_n=\frac{1}{2}\bigl(E_{n+1,n-1}+E_{n+2,n}\bigr),
\\
 h_0=-E_{1,1}-E_{2,2}+E_{2n-1,2n-1}+E_{2n,2n},\\
 h_i=E_{i,i}-E_{i+1,i+1}+E_{2n-i,2n-i}-E_{2n+1-i,2n+1-i},\quad 1\le i\le n-1,\\
  h_n=E_{n-1,n-1}+E_{n,n}-E_{n+1,n+1}-E_{n+2,n+2},
\\
h_\rho=\diag\bigl({-}n+1,\dots,0,0,\dots,n-1\bigr),
\quad
\ell=\diag\bigl(\mu_1,\dots,\mu_n,-\mu_n,\dots,-\mu_1\bigr).
\end{gather*}

\section{Proof of Theorem \ref{QQ-A_n}}\label{app2}

In this section we give a proof of Theorem \ref{QQ-A_n}.
Let $\mathfrak{g}=A^{(1)}_n$, and set
\begin{gather*}
\mathcal{A}_{i,z}=\mathcal{Y}_{i,q^{-1}z}^{-1}\mathcal{Y}_{i,qz}^{-1}\mathcal{Y}_{i-1,z}\mathcal{Y}_{i+1,z}, \qquad i=1,\ldots,n.
\end{gather*}
One can rewrite $\mathcal{Q}^{(1)}_{i,z}$ as
\begin{gather*}
\mathcal{Q}^{(1)}_{i,z} =\mathcal{Y}^{-1}_{n,q^{n+1}z}\sum_{k_{i+1},\dots, k_{n+1}}
\prod^{n}_{j=i}\left(\left(\frac{x_{j+1}}{x_i}\right)^{k_{j+1}}\mathcal{A}^{-1}_{j,q^{j-2\sum^{n+1}_{\mu=j+1}k_{\mu}}z}\right).
\end{gather*}

For $i\in\mathcal{J}$, we def\/ine
\begin{gather*}
 \Delta^{(0)}_{i,z}=\mathcal{Y}_{n,q^{n+1}z}\mathcal{Q}^{(1)}_{i,z},\\
 \Delta^{(a)}_{i,z}=\mathcal{A}_{n+1-a,q^{n+1-a}z}^{-1}\left(\Delta^{(a-1)}_{i,z}-\Delta^{(a-1)}_{i,q^{-2}z}\frac{x_{n+2-a}}{x_i}\right),
\qquad a=1,\ldots,n.
\end{gather*}

By the above def\/inition and the formula of $\mathcal{Q}^{(1)}_{i,z}$, we have
the following lemma.

\begin{Lemma}\label{simplify-A}
For $i\in\mathcal{J}$ and $0\leq a\leq n$, we have
\begin{gather*}
 \Delta^{(a)}_{i,z}=\begin{cases}
\sum\limits_{k_{i+1},\dots,k_{n+1-a}}\prod^{n-a}\limits_{j=i}\left(\left(\frac{x_{j+1}}{x_{i}}\right)^{k_{j+1}}
\mathcal{A}^{-1}_{j,q^{j-2\sum^{n+1-a}_{\mu=i}k_{\mu}}z}\right), &a< n+1-i,\\
1,&a=n+1-i,\\
0, &a>n+1-i.
\end{cases}
\end{gather*}
\end{Lemma}

\begin{Proposition}\label{New-Rep}
For an element $J=(j_1,j_2,\ldots,j_a)\in\mathcal{J}^a$, we have
\begin{gather*}
\mathcal{Q}^{(a)}_{J,z}
= \mathcal{Y}_{n+1-a,q^{n+1}z}^{-1}
\begin{vmatrix}
\Delta^{(0)}_{j_1,q^{1-a}z}x_{j_1}^0 &\Delta^{(0)}_{j_2,q^{1-a}z}x_{j_2}^{0}&\cdots&
\Delta^{(0)}_{j_a,q^{1-a}z}x_{j_a}^{0}\\
\Delta^{(1)}_{j_1,q^{3-a}z}x_{j_1}&\Delta^{(1)}_{j_2,q^{3-a}z}x_{j_2}&\cdots&\Delta^{(1)}_{j_a,q^{3-a}z}x_{j_a}^{}\\
\vdots&\vdots&&\vdots\\
\Delta^{(a-1)}_{j_1,q^{a-1}z}x_{j_1}^{a-1}&\Delta^{(a-1)}_{j_2,q^{a-1}z}x_{j_2}^{a-1}&
\cdots&\Delta^{(a-1)}_{j_a,q^{a-1}z}x_{j_a}^{a-1}
\end{vmatrix}.
\end{gather*}
\end{Proposition}

\begin{proof}
From the def\/inition we have
\begin{gather*}
\mathcal{Q}^{(a)}_{J,z}= \prod^a_{k=1}\mathcal{Y}_{n,q^{n-a+2k}z}^{-1}
\begin{vmatrix}
\Delta^{(0)}_{j_1,q^{1-a}z}x_{j_1}^0 &\Delta^{(0)}_{j_2,q^{1-a}z}x_{j_2}^{0}&\cdots&
\Delta^{(0)}_{j_a,q^{1-a}z}x_{j_a}^{0}\\
\Delta^{(0)}_{j_1,q^{3-a}z}x_{j_1}&\Delta^{(0)}_{j_2,q^{3-a}z}x_{j_2}&\cdots&\Delta^{(0)}_{j_a,q^{3-a}z}x_{j_a}^{}\\
\vdots&\vdots&&\vdots\\
\Delta^{(0)}_{j_1,q^{a-1}z}x_{j_1}^{a-1}&\Delta^{(0)}_{j_2,q^{a-1}z}x_{j_2}^{a-1}&
\cdots&\Delta^{(0)}_{j_a,q^{a-1}z}x_{j_a}^{a-1}
\end{vmatrix}.
\end{gather*}
For $k=2,\ldots,a$,
we subtract the $(k-1)$-th row multiplied by $x_{n+1}$
from the $k$-th row
and then extract the factor $\mathcal{A}^{-1}_{n,q^{n-a-1+2k}z}$. Thus we get
\begin{gather*}
\mathcal{Q}^{(a)}_{I,z}
= \prod^a_{k=1}\mathcal{Y}_{n,q^{n-a+2k}z}^{-1}\prod^{a}_{k=2}\mathcal{A}^{-1}_{n,q^{n-a-1+2k}z}\\
\hphantom{\mathcal{Q}^{(a)}_{I,z}=}{}
\times
\begin{vmatrix}
\Delta^{(0)}_{j_1,q^{1-a}z}x_{j_1}^0 &\Delta^{(0)}_{j_2,q^{1-a}z}x_{j_2}^{0}&\cdots&
\Delta^{(0)}_{j_a,q^{1-a}z}x_{j_a}^{0}\\
\Delta^{(1)}_{j_1,q^{3-a}z}x_{j_1}&\Delta^{(1)}_{j_2,q^{3-a}z}x_{j_2}&\cdots&\Delta^{(1)}_{j_a,q^{3-a}z}x_{j_a}^{}\\
\vdots&\vdots&&\vdots\\
\Delta^{(1)}_{j_1,q^{a-1}z}x_{j_1}^{a-1}&\Delta^{(1)}_{j_2,q^{a-1}z}x_{j_2}^{a-1}&
\cdots&\Delta^{(1)}_{j_a,q^{a-1}z}x_{j_a}^{a-1}
\end{vmatrix}.
\end{gather*}
Repeating the above steps in the order $b=2,\ldots,a-1$,
by replacing  $x_{n+1}$ and $\mathcal{A}^{-1}_{n,q^{n-a-1+2k}z}$
with $x_{n+2-b}$ and $\mathcal{A}^{-1}_{n+1-b,q^{n-a-b+2k}z}$
($k=b+1,\dots,a$), we get
\begin{gather*}
\mathcal{Q}^{(a)}_{I,z}
= \prod^a_{k=1}\mathcal{Y}_{n,q^{n-a+2k}z}\prod^{a-1}_{b=1}\prod^{a}_{k=b+1}\mathcal{A}^{-1}_{n+1-b,q^{n-a-b+2k}z}\\
\hphantom{\mathcal{Q}^{(a)}_{I,z}=}{}
\times
\begin{vmatrix}
\Delta^{(0)}_{j_1,q^{1-a}z}x_{j_1}^0 &\Delta^{(0)}_{j_2,q^{1-a}z}x_{j_2}^{0}&\cdots&
\Delta^{(0)}_{j_a,q^{1-a}z}x_{j_a}^{0}\\
\Delta^{(1)}_{j_1,q^{3-a}z}x_{j_1}&\Delta^{(1)}_{j_2,q^{3-a}z}x_{j_2}&\cdots&\Delta^{(1)}_{j_a,q^{3-a}z}x_{j_a}\\
\vdots&\vdots&&\vdots\\
\Delta^{(a-1)}_{j_1,q^{a-1}z}x_{j_1}^{a-1}&\Delta^{(a-1)}_{j_2,q^{a-1}z}x_{j_2}^{a-1}&
\cdots&\Delta^{(a-1)}_{j_a,q^{a-1}z}x_{j_a}^{a-1}
\end{vmatrix}.
\end{gather*}
The proof is over by noting that
\begin{gather*}
\prod^a_{k=1}\mathcal{Y}_{n,q^{n-a+2k}z}^{-1}\prod^{a-1}_{b=1}\prod^{a}_{k=b+1}\mathcal{A}^{-1}_{n+1-b,q^{n-a-b+2k}z}=\mathcal{Y}_{n+1-a,q^{n+1}z}^{-1}.
\tag*{\qed}
\end{gather*}
\renewcommand{\qed}{}
\end{proof}

Now we give a proof of Theorem \ref{QQ-A_n}.

\begin{proof}[Proof of Theorem \ref{QQ-A_n}] By Proposition \ref{New-Rep}, for $J=(1,2,\ldots,n+1)$, we have
\begin{gather*}
\mathcal{Q}^{(n+1)}_{J,z}
=
\begin{vmatrix}
\Delta^{(0)}_{1,q^{-n}z}x_1^0 &\Delta^{(0)}_{2,q^{-n}z}x_2^{0}&\cdots&
\Delta^{(0)}_{n+1,q^{-n}z}x_{n+1}^{0}\\
\Delta^{(1)}_{1,q^{2-n}z}x_1&\Delta^{(1)}_{2,q^{2-n}z}x_2&\cdots&\Delta^{(1)}_{n+1,q^{2-n}z}x_{n+1}\\
\vdots&\vdots&&\vdots\\
\Delta^{(n)}_{1,q^{n}z}x_1^{n}&\Delta^{(n)}_{2,q^{n}z}x_2^n&
\cdots&\Delta^{(n)}_{n+1,q^{n}z}x_{n+1}^n
\end{vmatrix}.
\end{gather*}
By Lemma \ref{simplify-A}, we get
\begin{gather*}
\mathcal{Q}^{(n+1)}_{J,z}
= (-1)^{(n+1)n/2}\prod^{n}_{b=1}x_b^{n+1-b},
\end{gather*}
which implies the result.
\end{proof}

\section[Proofs of identities for characters in the case $B^{(1)}_n$]{Proofs of identities for characters in the case $\boldsymbol{B^{(1)}_n}$}\label{app3}

In this section, we
set $\mathfrak{g}=B^{(1)}_n$.
First we introduce some notation.
For $i,j\in\mathcal{J}$ and an element $J\in\mathcal{J}^a$, we def\/ine
\begin{gather*}
 h_{J}=\prod_{k,\bar{k}\in J,\atop k\prec\bar{k}}\big(x_k^{\frac{1}{2}}+x_k^{-\frac{1}{2}}\big),
\qquad \tilde{x}_{i,J}=
\prod_{k\in J,\atop k\prec i}\big({-}x_i^{\frac{1}{2}\delta_{i,\bar{k}}-1}\big)
\prod_{k\in J^*,\atop i\prec k}\big(x_k^{1-\frac{1}{2}\delta_{k,\bar{i}}}\big), \\
 g_{i,j}=f_{i,j}x_j^{1-\frac{1}{2}\delta_{j,\bar{i}}}=-g_{j,i},
\qquad g_{i,J}=\prod_{k\in J,\atop k\neq i}g_{i,k},
\end{gather*}
where $f_{i,j}$ is def\/ined in Subsection~\ref{section4.2}.

By \eqref{Character-Formula} and the above notation, one has
\begin{gather*}
\bar{Q}^{(a)}_{J} =\prod_{i\in J,j\in J^*,\atop i\prec j}f_{i,j}
\prod_{i,j\in J,\atop i\prec j}\left((x_j-x_i)f_{i,j}\right)
=X^{-1}_{J}h_{J}\prod_{i\in J,j\in J^*,\atop i\prec j}f_{i,j},
\end{gather*}
where $X_J$ is given in~\eqref{notation-B_n}.
A direct computation gives
\begin{gather*}
 \prod_{i\in\overline{J^*},j\in\overline{J},\atop i\prec j}f_{i,j}=\prod_{i\in J,j\in J^*,\atop i\prec j}f_{i,j},\qquad
 \prod_{j\in\overline{J},\bar{i}\prec j}f_{\bar{i},j}=\prod_{j\in J,j\prec i}f_{j,i}
,\\
 \frac{\prod\limits_{i\in J-j_k,j\in (J-j_k)^*,\atop i\prec j}f_{i,j}}{\prod\limits_{i\in J,j\in J^*,\atop i\prec j}f_{i,j}}=\frac{\prod\limits_{i\in J,\atop i\prec j_k}g_{j_k,j}}{\prod\limits_{j\in J^*,\atop j_k\prec j}g_{j_k,j}}\tilde{x}_{j_k,J}.
\end{gather*}
Furthermore, for two elements $J_1$ and $J_2$ of $\mathcal{J}^n$ and $j_k\in \underset{\sim}{J_2}\cap\underset{\sim}{J^*_1}$, we have
\begin{gather*}
\frac{x_{j_k,J_1}\tilde{x}_{j_k,J_2}}{x_{j_k,J_2}\tilde{x}_{j_k,J_1}}=-x^2_{j_k}x_{J_1}x^{-1}_{J_2}.
\end{gather*}
Using the above formulas, the specialization of
Conjecture~\ref{conj-B_n} to the characters $\bar{\mathcal{Q}}^{(a)}_J$
reduces to the following proposition.

\begin{Proposition}
For two increasing elements
$J_1=(i_1,\ldots,i_n), \, J_2=(j_1,\ldots,j_n)\in \mathcal{J}^n$, we have
\begin{gather}
 x_{J_1}^{-\frac{1}{2}}h_{\overline{J^*_1}}h_{J_2}
-x_{J_2}^{-\frac{1}{2}}h_{\overline{J^*_2}}h_{J_1}=x_{J_1}^{-1}x_{J_2}^{-1}\sum_{j\in J_2\cap J_1^*}x_{J_1}^{\frac{3}{2}}x_{j}^{\frac{3}{2}}
\frac{g_{j,J_1^*}}{g_{j,J_2^*}}h_{J_2-j}h_{\overline{J_1^*}-\overline{j}}.\label{propAppenC}
\end{gather}
\end{Proposition}

\begin{proof} Generally let
\begin{gather*}
 T_1=\left\{j\in\underset{\sim}{J_1}\cap\underset{\sim}{J_2^*}\left|\, \bar{j}\notin J_1,\bar{j}\in J_2^*\;\text{or}\;\bar{j}\in\underset{\sim}{J_1},\bar{j}\notin\underset{\sim}{J_2^*}\right.\right\},\\
 T_2=\left\{j\in\underset{\sim}{J_1}\cap\underset{\sim}{J_2^*}\left| \bar{j}\notin\underset{\sim}{J_1},\bar{j}\notin\underset{\sim}{J_2^*}\right.\right\},\qquad
T_3=\left\{j\in\underset{\sim}{J_1}\cap\underset{\sim}{J_2^*}\left| \bar{j}\in\underset{\sim}{J_1}\cap\underset{\sim}{J_2^*}, j\prec\bar{j}\right.\right\},\\
 S_1=\left\{j\in\underset{\sim}{J_2}\cap\underset{\sim}{J_1^*}\left| \bar{j}\notin\underset{\sim}{J_2},\bar{j}\in\underset{\sim}{J_1^*},\text{or}\;
\bar{j}\in\underset{\sim}{J_2},\bar{j}\notin\underset{\sim}{J_1^*}\right.\right\},\\
 S_2=\left\{j\in\underset{\sim}{J_2}\cap\underset{\sim}{J_1^*}\left| \bar{j}\notin\underset{\sim}{J_2},\bar{j}\notin\underset{\sim}{J_1^*}\right.\right\},\qquad
S_3=\left\{j\in\underset{\sim}{J_2}\cap\underset{\sim}{J_1^*}\left| \bar{j}\in\underset{\sim}{J_2}\cap\underset{\sim}{J_1^*}, j\prec\bar{j}\right.\right\},
\end{gather*}
then we have $\underset{\sim}{J_1}\cap\underset{\sim}{J_2^*}=\bigcup\limits_{1\leq i\leq 3}T_i\cup\overline{T}_3$, $\underset{\sim}{J_2}\cap\underset{\sim}{J_1^*}=\bigcup\limits_{1\leq i\leq 3}S_i\cup\overline{S}_3$,
and $\sharp (\underset{\sim}{J_1}\cap\underset{\sim}{J_2^*})
=\sharp (\underset{\sim}{J_2}\cap\underset{\sim}{J_1^*})$. Here for a subset $S\subseteq\mathcal{J}$, we set $\overline{S}=\{\bar{i}\mid i\in S\}$.

Let
\begin{gather*}
F_{J_1,J_2}(x_j)=x_{j}^{\frac{3}{2}}x_{J_1}x^{-1}_{J_2}\frac{h_{J_2-j}h_{\bar{J}_1^*-\bar{j}}}{h_{\bar{J}^*_1}h_{J_2}}
\frac{g_{j,J_1^*}}{g_{j,J_2^*}}.
\end{gather*}
We also introduce a function
\begin{gather*}
f(z)=\frac{1}{(1+z)}\frac{\prod\limits_{k\in T_1}(z-x_k)}{\prod\limits_{k\in S_1}(z-x_k)}\prod_{k\in S_2}\frac{zx_k-1}{z-x_k}\frac{\prod\limits_{k\in T_3}(z-x_k)(zx_k-1)}{\prod\limits_{k\in S_3}(z-x_k)(zx_k-1)},
\end{gather*}
which has poles only at
$z=-1,\infty$ and $x_j$
($j\in \underset{\sim}{J_2}\cap\underset{\sim}{J_1^*}=\bigcup\limits_{1\leq i\leq 3}S_i\cup\overline{S}_3$).

A direct computation gives, for
$j\in S_3$
\begin{gather*}
 \text{Res}\big(f,x^{-1}_j\big)
= \frac{1}{x_j(1+x_j^2)(1-x_j)}\frac{\prod\limits_{k\in T_1}(1-x_kx_j)}{\prod\limits_{k\in S_1}(1-x_kx_j)}\prod_{k\in S_2}\frac{x_k-x_j}{1-x_kx_j}\frac{\prod\limits_{k\in T_3}(x_j-x_k)(x_jx_k-1)}{\prod\limits_{k\in S_3,\atop k\neq j}(x_j-x_k)(x_jx_k-1)}\\
\hphantom{\text{Res}\big(f,x^{-1}_j\big)}{} = cF_{J_1,J_2}(x^{-1}_j),
\end{gather*}
where $c=\frac{x_{T_3}x_{S_2}}{x_{S_3}}$.
Similarly we obtain
\begin{gather*}
\text{Res}(f,x_j)=cF_{J_1,J_2}(x_j),\qquad  j\in S_i, \quad i=1,2,3.
\end{gather*}
On the other hand,
by setting
$h_i=x_i^{\frac{1}{2}}+x_i^{-\frac{1}{2}}$, it is easy to see that
\begin{gather*}
 \text{Res}\; (f,-1)=\frac{x_{T_3}}{x_{S_3}}\left(\frac{x_{T_1}}{x_{S_1}}\right)^{\frac{1}{2}}\frac{\prod\limits_{k\in T_1}h_k}{\prod\limits_{k\in S_1}h_k}\frac{\prod\limits_{k\in T_3}h_k^2}{\prod\limits_{k\in S_3}h_k^2}, \qquad\text{Res}\;(f,\infty)=-\frac{x_{T_3}x_{S_2}}{x_{S_3}}=-c.
\end{gather*}
Since the right-hand side of \eqref{propAppenC} equals
$h_{\bar{J}^c_1}h_{J_2}x_{J_1}^{-\frac{1}{2}}\sum\limits_{j\in\underset{\sim}{J_2}\cap\underset{\sim}{J_1^*}}F_{J_1,J_2}(x_j)$,
it simplif\/ies to
\begin{gather*}
 h_{\bar{J}^*_1}
h_{J_2}x_{J_1}^{-\frac{1}{2}}
(-c^{-1})\left(-c
+\frac{x_{T_3}}{x_{S_3}}\left(\frac{x_{T_1}}{x_{S_1}}\right)^{\frac{1}{2}}\frac{\prod\limits_{k\in T_1}h_k}{\prod\limits_{k\in S_1}h_k}\frac{\prod\limits_{k\in T_3}h_k^2}{\prod\limits_{k\in S_3}h_k^2}\right).
\end{gather*}
Now it is suf\/f\/icient to prove
\begin{gather*}
 \left(\frac{x_{T_1}}{x_{S_1}}\right)^{\frac{1}{2}}\frac{\prod\limits_{k\in T_1}h_k}{\prod\limits_{k\in S_1}h_k}\frac{\prod\limits_{k\in T_3}h_k^2}{\prod\limits_{k\in S_3}h_k^2}=x_{S_2}\frac{x_{J_1}^{\frac{1}{2}}}{x_{J_2}^{\frac{1}{2}}}\frac{h_{\bar{J}^*_2}h_{J_1}}{h_{\bar{J}^*_1}h_{J_2}}
\end{gather*}
which is immediate. We f\/inish the proof.
\end{proof}

\section[Proofs for identities of characters of case $C^{(1)}_n$ and $D^{(1)}_n$]{Proofs for identities of characters of case $\boldsymbol{C^{(1)}_n}$ and $\boldsymbol{D^{(1)}_n}$}\label{app4}

In this section we consider $\mathfrak{g}=C^{(1)}_n, D^{(1)}_n$.
By \eqref{Character-Formula} and Conjecture~\ref{spin-character-C_n},
Conjectures \ref{conj-C_n-1}--\ref{conj-C_n-4} reduce to the following identities
\begin{gather}
 (-1)^{\frac{u(u+1)}{2}}x_{\pmb{\sigma}}^2\prod_{1\leq k\leq u,\atop 1\leq l\leq u+1}
\frac{1}{(1-x_{\sigma_k}x_{\eta_l})(x_{\sigma_k}-x_{\eta_l})}\prod_{1\leq k\leq u+1}(1-x_{\eta_k}^2)
\nonumber\\
\qquad{} = \sum(-1)^{\frac{s(s+1)}{2}}x_{\pmb{\gamma}}^2\prod_{1\leq k<l\leq s}\frac{1}{1-x_{\gamma_k}x_{\gamma_l}}
\prod_{1\leq k<l\leq t}\frac{1}{1-x_{\delta_k}x_{\delta_l}}
\prod_{1\leq k\leq s,\atop 1\leq l\leq t}\frac{1}{x_{\gamma_k}-x_{\delta_l}},
\label{chC1}\\
 (-1)^{\frac{u(u-1)}{2}}x_{\pmb{\sigma}}^2\prod_{1\leq k\leq u,\atop 1\leq l\leq u+2}
\frac{1}{(1-x_{\sigma_k}x_{\eta_l})(x_{\sigma_k}-x_{\eta_l})}
\prod_{1\leq k\leq u}(1-x_{\sigma_k})\prod_{1\leq k\leq u+2}(1+x_{\eta_k})
\nonumber\\
\qquad{} = \sum(-1)^{\frac{s(s-1)}{2}}x_{\pmb{\gamma}}^2\prod_{1\leq k<l\leq s}
\frac{1}{1-x_{\gamma_k}x_{\gamma_l}}\prod_{1\leq k<l\leq t}\frac{1}{1-x_{\delta_k}x_{\delta_l}}
\prod_{1\leq k\leq s,\atop 1\leq l\leq t}\frac{1}{x_{\gamma_k}-x_{\delta_l}},
\label{chC2}\\
 (-1)^{\frac{u(u+1)}{2}}x_{\pmb{\sigma}}\prod_{1\leq k\leq u,\atop 1\leq l\leq u}
\frac{1}{(1-x_{\sigma_k}x_{\eta_l})(x_{\sigma_k}-x_{\eta_l})}\prod_{1\leq k\leq u}(1-x^2_{\eta_k})
\nonumber\\
\qquad{} = \sum(-1)^{\frac{s(s+1)}{2}}x_{\pmb{\gamma}}\prod_{1\leq k<l\leq s}\frac{1}{1-x_{\gamma_k}x_{\gamma_l}}
\prod_{1\leq k<l\leq t}\frac{1}{1-x_{\delta_k}x_{\delta_l}}
\prod_{1\leq k\leq s,\atop 1\leq l\leq t}\frac{1}{x_{\gamma_k}-x_{\delta_l}},
\label{chC3}\\
 (-1)^{\frac{u(u-1)}{2}}x_{\pmb{\sigma}}\prod_{1\leq k\leq u,\atop 1\leq l\leq u+1}
\frac{1}{(1-x_{\sigma_k}x_{\eta_l})(x_{\sigma_k}-x_{\eta_l})}
\prod_{1\leq k\leq u}(1-x_{\sigma_k})\prod_{1\leq k\leq u+1}(1+x_{\eta_k})
\nonumber\\
\qquad{} = \sum(-1)^{\frac{s(s-1)}{2}}x_{\pmb{\gamma}}\prod_{1\leq k<l\leq s}
\frac{1}{1-x_{\gamma_k}x_{\gamma_l}}\prod_{1\leq k<l\leq t}\frac{1}{1-x_{\delta_k}x_{\delta_l}}
\prod_{1\leq k\leq s,\atop 1\leq l\leq t}\frac{1}{x_{\gamma_k}-x_{\delta_l}}.
\label{chC4}
\end{gather}
In these formulas, $\pmb{\sigma}$, $\pmb{\eta}$ are disjoint subsets of $\{1,\ldots,n\}$,
such that $(\sharp\,\pmb{\sigma},\sharp\,\pmb{\eta})=(u,u+1), (u,u+2)$, $(u,u),(u,u+1)$ respectively.
The sum is taken over partitions of
$\pmb{\sigma}\cup\pmb{\eta}$ into subsets $\pmb{\gamma}$, $\pmb{\delta}$
satisfying $\pmb{\sigma}\subseteq\pmb{\gamma}$,
where we set $s=\sharp\,\pmb{\gamma}$, $t=\sharp\,\pmb{\delta}$.

Similarly, in the case of $\mathfrak{g}=D^{(1)}_n$,
by \eqref{Character-Formula} and Conjecture~\ref{spin-character-D_n}
the relevant identities read as follows.
\begin{gather}
 (-1)^{\frac{u(u+1)}{2}}x_{\pmb{\sigma}}\prod_{1\leq k\leq u,\atop 1\leq l\leq u+2}
\frac{1}{(1-x_{\sigma_k}x_{\eta_l})(x_{\sigma_k}-x_{\eta_l})}\prod_{1\leq k\leq u}\big(1-x^2_{\sigma_k}\big)
\nonumber
\\
\qquad{} = \sum(-1)^{\frac{t(t+1)}{2}}x_{\pmb{\delta}}\prod_{1\leq k<l\leq s}
\frac{1}{1-x_{\gamma_k}x_{\gamma_l}}\prod_{1\leq k<l\leq t}\frac{1}{1-x_{\delta_k}x_{\delta_l}}
\prod_{1\leq k\leq s,\atop 1\leq l\leq t}\frac{1}{x_{\gamma_k}-x_{\delta_l}},
\label{chD1}\\
 (-1)^{\frac{u(u-1)}{2}}\prod_{1\leq k\leq u,\atop 1\leq l\leq u+2}
\frac{1}{(1-x_{\sigma_k}x_{\eta_l})(x_{\sigma_k}-x_{\eta_l})}\prod_{1\leq k\leq u}\big(1-x^2_{\sigma_k}\big)
\nonumber
\\
\qquad{}= \sum(-1)^{\frac{s(s-1)}{2}}\prod_{1\leq k<l\leq s}\frac{1}{1-x_{\gamma_k}x_{\gamma_l}}
\prod_{1\leq k<l\leq t}\frac{1}{1-x_{\delta_k}x_{\delta_l}}
\prod_{1\leq k\leq s,\atop 1\leq l\leq t}\frac{1}{x_{\gamma_k}-x_{\delta_l}}.
\label{chD2}
\end{gather}
Here  $(\sharp\,\pmb{\sigma},\sharp\,\pmb{\eta})=(u,u+2),(u,u+1)$ respectively,
 $s=\sharp\,\pmb{\gamma}$, $t=\sharp\,\pmb{\delta}$, and the sum is taken over the partitions
satisfying $\pmb{\sigma}\subseteq\pmb{\gamma}$ keeping f\/ixed the parity of $s$.

These formulas are related by specialization of variables:
\eqref{chC1} is obtained from~\eqref{chC3} by setting $x_{\sigma_u}=0$,
\eqref{chC2} is obtained from~\eqref{chC4} by setting $x_{\sigma_u}=0$,
and
\eqref{chC4} is obtained from~\eqref{chC3} by setting $x_{\sigma_u}=1$.
Likewise~\eqref{chD1} is obtained from~\eqref{chD2} by setting $x_{\sigma_u}=\infty$.
Hence it is suf\/f\/icient to deal only with~\eqref{chC3} and~\eqref{chD2}.

Def\/ine polynomials in $x=(x_1,\ldots,x_m)$ and $y=(y_1,\ldots,y_n)$,
\begin{gather*}
F^{(m,n)}(x,y)
 =\prod\limits_{1\le i<j\le m}(1-x_ix_j)\prod\limits_{1\le i<j\le n}(1-y_iy_j)(y_j-y_i),
\\
G^{(m,n)}_{J_1,J_2}(x,y) =
\prod\limits_{1\le i\le m, \atop j\in J_1}(1-x_iy_j)\prod\limits_{1\le i\le m,\atop j\in J_2}(x_i-y_j)
\\
\hphantom{G^{(m,n)}_{J_1,J_2}(x,y) =}{}  \times
\prod\limits_{i,j\in J_1,\atop i<j}(y_j-y_i)
\prod\limits_{i\in J_1, j\in J_2}(1-y_iy_j)
\prod\limits_{i,j\in J_2, \atop i<j}(y_i-y_j).
\end{gather*}
Set also $\delta(J_1,J_2)=\sharp\,\{(i,j)\in J_1\times J_2\mid i<j\}$.
Then  \eqref{chC3} and \eqref{chD2} can be rewritten respectively as follows:
\begin{gather}
 \prod_{j=1}^n\big(1-y_j^2\big)\cdot F^{(n,n)}(x,y)
=\sum_{J_1,J_2}(-1)^{\delta(J_1,J_2)}\prod_{j\in J_2}y_j\cdot G^{(n,n)}_{J_1,J_2}(x,y),
\label{cd1}
\\
 (1+\epsilon)\prod_{i=1}^{n-1}\big(1-x_i^2\big)\cdot F^{(n-1,n)}(x,y)
=\sum_{J_1,J_2}\epsilon^{\sharp J_2}
(-1)^{\delta(J_1,J_2)}G^{(n-1,n)}_{J_1,J_2}(x,y).
\label{cd2}
\end{gather}
Here the sum in the right-hand sides are taken over all partitions $J_1\sqcup J_2$ of $\{1,\ldots,n\}$,
and $\epsilon=\pm1$.
For the second, we have added/subtracted the original sums which have
a parity restriction.

In order to show \eqref{cd1}, we prove it in a slightly more general form.

\begin{Proposition}\label{Did1}
Suppose $n\ge 2$, $0\le m\le n$ and $m\equiv n\pmod 2$.
  Then
\begin{gather}
 \delta_{m,n}\prod_{j=1}^n\big(1-y_j^2\big)\cdot F^{(m,n)}(x,y)
=\sum_{J_1,J_2}(-1)^{\delta(J_1,J_2)}\prod_{j\in J_2}y_j\cdot G^{(m,n)}_{J_1,J_2}(x,y).
\label{FG1}
\end{gather}
\end{Proposition}

\begin{proof}
Denote by $R^{(m,n)}$ the right-hand side of \eqref{FG1}, and set
\begin{gather*}
L^{(m,n)}=\prod_{j=1}^n\big(1-y_j^2\big)\cdot F^{(m,n)}(x,y),\qquad
h^{(m)}(x,z)=\prod_{i=1}^m(1-x_iz).
\end{gather*}
Let us rewrite $R^{(m,n)}$ as follows:
\begin{gather*}
 \prod_{j=1}^ny_j^{-(m+n)/2}\cdot R^{(m,n)}
 =\sum_{J_1,J_2}\!(-1)^{\sharp J_2}
\prod_{i\in J_1}\big(y_i^{-(m+n)/2}h^{(m)}(x,y_i)\big)\!\prod_{j\in J_2}\big(y_j^{(m+n)/2}h^{(m)}(x,y^{-1}_j)\big)
\\
\qquad{}  \times \prod_{i,j\in J_1,\atop i<j}(y_j-y_i)
 \prod_{i\in J_1,j\in J_2,\atop i<j}\big(y_j^{-1}-y_i\big)
 \prod_{i\in J_1,j\in J_2,\atop i>j}\big(y_i-y_j^{-1}\big)
\prod_{i,j\in J_2,\atop i<j}\big(y_j^{-1}-y_i^{-1}\big)
\\
\qquad{}   =\sum_{\varepsilon_1,\ldots,\varepsilon_n=\pm1}
\prod_{j=1}^n\left(
\varepsilon_j y_j^{-\varepsilon_j(m+n)/2}h^{(m)}(x,y_j^{\varepsilon_j})\right)
\prod_{1\le i<j\le n}(y_j^{\varepsilon_j}-y_i^{\varepsilon_i})\\
\qquad{}
  =\det\left(y_i^{-(m+n)/2+j-1}h^{(m)}(x,y_i)
-y_i^{(m+n)/2-j+1}h^{(m)}\big(x,y_i^{-1}\big)\right)_{1\le i,j\le n}.
\end{gather*}
It is easy to see that the last expression is $0$ if $m=0$ or $m=1$.
Assume $m\ge2$.
We prove~\eqref{FG1} by induction on $m$.
From the skew-symmetry of the determinant, $R^{(m,n)}$
is divisible by $\prod_{j=1}^n(1-y_j^2)\prod_{1\le i<j\le n}(1-y_iy_j)(y_j-y_i)$.
Let us show that it is divisible also by $\prod_{1\le i<j\le m}(1-x_ix_j)$.

If we set $x_{m}=x_{m-1}^{-1}$, then $h^{(m)}(x,z)=z\big(z+z^{-1}-x_{m-1}-x_{m-1}^{-1}\big)h^{(m-2)}(x'',z)$, where
$x''=(x_1,\ldots,x_{m-2})$. Hence $R^{(m,n)}$ becomes proportional to
\begin{gather*}
\det\left(y_i^{-(m-2+n)/2+j-1}h^{(m-2)}(x'',y_i)-y_i^{(m-2+n)/2-j+1}h^{(m-2)}\big(x'',y_i^{-1}\big)\right)_{1\le i,j\le n},
\end{gather*}
which is $0$ by the induction hypothesis. This shows that $R^{(m,n)}$ is divisible by $L^{(m,n)}$.

The total degree of $L^{(m,n)}$ is $l^{(m,n)}=m(m-1)+n(3n+1)/2$, while that of $R^{(m,n)}$ is at most
$r^{(m,n)}=n^2/2+mn+(m+n)^2/4-m/2$.
Notice that $l^{(m,n)}>r^{(m,n)}$ for $n>m$ and $l^{(n,n)}=r^{(n,n)}$.
From this we conclude that if $n>m$ then
 $R^{(m,n)}=0$, and if $n=m$ then $R^{(n,n)}$ is a constant multiple of $L^{(n,n)}$.
The constant is shown to be $1$ by setting $x=0$ and using the Weyl denominator formula of type $C_n$,
\begin{gather*}
\prod_{j=1}^ny_j^n\cdot
\det\left(y_i^{-n+j-1}-y_i^{n-j+1}\right)_{1\le i,j\le n}=
\prod_{j=1}^n\big(1-y_j^2\big)\prod_{1\le i<j\le n}(1-y_iy_j)(y_j-y_i).
\tag*{\qed}
\end{gather*}
\renewcommand{\qed}{}
\end{proof}

\begin{Proposition}
Suppose $n\ge 3$ and $0\le m\le n-1$. Then for $\epsilon=\pm1$ we have
\begin{gather}
 (1+\epsilon)\delta_{m,n-1}\prod_{i=1}^{m}\big(1-x_i^2\big)\cdot F^{(m,n)}(x,y)
=\sum_{J_1,J_2}\epsilon^{\sharp J_2}
(-1)^{\delta(J_1,J_2)}G^{(m,n)}_{J_1,J_2}(x,y) .
\label{FG2}
\end{gather}
\end{Proposition}

\begin{proof}
We denote the right hand-side of \eqref{FG2} by $R^{(m,n)}_\epsilon$, and set
\begin{gather*}
K^{(m,n)}=\prod_{i=1}^{m}(1-x_i^2)\cdot F^{(m,n)}(x,y).
\end{gather*}
By a calculation similar to the one in the proof of Proposition \ref{Did1}
we f\/ind
\begin{gather*}
 \prod_{j=1}^ny_j^{-(m+n-1)/2}\cdot R^{(m,n)}_\epsilon
\\
 \qquad{}=\det\left(y_i^{-(m+n+1)/2+j}h^{(m)}(x,y_i)
+\epsilon' y_i^{(m+n+1)/2-j}h^{(m)}(x,y_i^{-1})\right)_{1\le i,j\le n},
\end{gather*}
where $\epsilon'=(-1)^{m+n-1}\epsilon $.
From this, it is easy to see that $R^{(0,n)}_\epsilon=0$.
As in Proposition \ref{Did1}, we prove the assertion by induction on $m$.
As before, one can verify that the left-hand side is divisible by
$\prod_{1\le i<j\le m}(1-x_ix_j)\prod_{1\le i<j\le n}(1-y_iy_j)(y_j-y_i)$.
Let us show that it is divisible also by $\prod_{i=1}^m(1-x_i^2)$.
If we set $x_m=\pm 1$, then we have
\begin{align*}
h^{(m)}(x,z)\Bigl|_{x_m=\pm 1}=\big(z^{-1/2}\mp z^{1/2}\big) \cdot z^{1/2}h^{(m-1)}(x',z),\qquad
x'=(x_1,\ldots,x_{m-1}).
\end{align*}
Hence the determinant becomes proportional to \begin{align*}
&
\det\left(y_i^{-(m+n)/2+j}h^{(m-1)}(x',y_i)\mp \epsilon'
y_i^{(m+n)/2-j}h^{(m-1)}\big(x,y_i^{-1}\big)\right)_{1\le i,j\le n},
\end{align*}
which is $0$ by the induction hypothesis.

The total degree of $K^{(m,n)}$ is $k^{(m,n)}=m(m+1)+3n(n-1)/2$, while
that of $R^{(m,n)}_\epsilon$ is at most $r^{(m,n)}_\epsilon=mn+n(n-1)/2+(m+n)^2/4-\delta/4$,
where $\delta=0,1$ is determined by $\delta\equiv m-n\pmod2$.
Then $k^{(m,n)}>r^{(m,n)}_\epsilon$ if $m\le n-2$, and $k^{(n-1,n)}=r^{(n-1,n)}_\epsilon$.
Hence we f\/ind that $R^{(m,n)}_\epsilon=0$ if $m\le n-2$, and
$R^{(n-1,n)}_\epsilon$ is a constant multiple of $K^{(n-1,n)}$.
The constant can be found by setting $x=0$ and
using the Weyl denominator formula of type $D_n$,
\begin{align*}
&\prod_{j=1}^ny_j^n\cdot
\det\left(y_i^{-n+j}+\epsilon\, y_i^{n-j}\right)_{1\le i,j\le n}=
(1+\epsilon) \prod_{1\le i<j\le n}(1-y_iy_j)(y_j-y_i).
\tag*{\qed}
\end{align*}
\renewcommand{\qed}{}
\end{proof}

\subsection*{Acknowledgements}

The author would like to express her deep gratitude
to Professor Michio Jimbo for his valuable
advices, great helps and kindness, and to Professor Hidetaka Sakai for his constant encouragement during the course of the present work.
She wishes to thank also Professor Junji Suzuki for helpful discussions. Finally the author would like to thank the anonymous referees for comments and suggestions that improved the paper greatly.

\pdfbookmark[1]{References}{ref}
\LastPageEnding

\end{document}